\documentclass[final]{siamltex}

\usepackage{amssymb,latexsym,amsmath}
\usepackage{epsfig}
\usepackage{caption}
\usepackage{caption}
 \usepackage{color}
 
\newcommand{\boldGamma}{{\bf \Gamma}}

\newcommand{\boldu}{{\bf u}}

\newcommand{\boldy}{{\bf y}}

\newcommand{\undergamma}{{\underline{\gamma}}}
\newcommand{\IN}{{\it N}}

\def\IR{I \kern-0.35em R\,}

\def\l1{L_1}
\def\l2{L_2}

\newlength{\noteWidth}
\long\def\notes#1{\ifinner
           {\tiny #1}
           \else
           \marginpar{\parbox[t]{\noteWidth}{\raggedright\tiny #1}}
       \fi\typeout{#1}}

       \newtheorem{thm}{\bf{Theorem}}[section]
       
       \newtheorem{lem}{\bf{Lemma}}[section]
       \newtheorem{fact}{\bf{Fact}}[section]
       
       \newtheorem{defn}{\bf{Definition}}[section]
  \newtheorem{remark}{Remark}
  
  \newtheorem{example}{Example}

\def\qed{\hfill $\diamond$}
\fussy

\begin{document}

\title{Convex Analysis in Decentralized Stochastic Control, Strategic Measures and Optimal Solutions}

\author{Serdar Y\"uksel and Naci Saldi
\thanks{Serdar Y\"uksel is with the Department of Mathematics and
    Statistics, Queen's University, Kingston, Ontario, Canada; email: yuksel@mast.queensu.ca. Naci Saldi is with the Coordinated Science Laboratory, University of Illinois at Urbana-Champaign, Illinois, USA; email: nsaldi@illinois.edu. This research was partially supported by the Natural Sciences and Engineering Research Council of Canada (NSERC). Part of the results of this paper is to appear at the IEEE Conference on Decision and Control (2016). }
}

\maketitle

\begin{abstract}
This paper is concerned with the properties of the sets of strategic measures induced by admissible team policies in decentralized stochastic control and the convexity properties in dynamic team problems. To facilitate a convex analytical approach, strategic measures for team problems are introduced. Properties such as convexity, compactness and Borel measurability under weak convergence topology are studied, and sufficient conditions for each of these properties are presented. These lead to existence of and structural results for optimal policies. It will be shown that the set of strategic measures for teams which are not classical is in general non-convex, but the extreme points of a relaxed set consist of deterministic team policies, which lead to their optimality for a given team problem under an expected cost criterion. Externally provided independent common randomness for static teams or private randomness for dynamic teams do not improve the team performance. The problem of when a sequential team problem is convex is studied and necessary and sufficient conditions for problems which include teams with a non-classical information structure are presented. Implications of this analysis in identifying probability and information structure dependent convexity properties are presented.
\end{abstract}

\begin{keywords}
Stochastic control, decentralized control, optimal control, convex analysis.
\end{keywords}

\begin{AMS}
93E03, 90B99, 49J55 	
\end{AMS}


\section{Introduction}

Team decision theory has its roots in control theory and economics. Marschak \cite{mar55} was perhaps the first to introduce the basic elements of teams, and to  provide the first steps toward the development of a {\it team theory}. Radner \cite{rad62} provided foundational results for static teams, establishing connections between person-by-person optimality, stationarity, and team-optimality \cite{marrad72}. Contributions of Witsenhausen \cite{wit71,wit75,wit88,WitsenStandard,WitsenhausenSIAM71} on dynamic teams and characterization of information structures have been crucial in the progress of our understanding of dynamic teams. We refer the reader to Section \ref{witsenInfoStructureReview}, where Witsenhausen's {\it intrinsic model}, and characterization of information structures are discussed in detail. Further discussion on design of information structures in the context of team theory is available in \cite{Arrow1985,vanZandt,YukselBasarBook}. 


Convexity is a very important property for optimization problems. A property related to convex analysis that is relevant in team problems is the characterization of the sets of strategic measures; these are the probability measures induced on the exogenous variables, and measurement and action spaces by admissible control policies. In the context of single decision maker control problems, such measures have been studied extensively in \cite{Schal,piunovskii1998controlled,dynkin1979controlled,feinberg1996measurability}. A study of strategic measures for team problems has not been made to our knowledge, and it will be observed in this paper that many of the properties that are natural for fully-observed single-decision-maker stochastic control problems, such as convexity, do not generally extend to a large class of stochastic team problems. On the other hand, new results on the existence and the structure of optimal team policies will be established through such a convex analytical formulation. 

Another dimension of a convex analytical approach for team problems is on the characterization of convexity of cost functions on the set of admissible team policies: In the context of decentralized control problems, under convexity, for global optimality of team policies it may suffice to search for person-by-person optimal solutions (through \cite{rad62},\cite{KraSpeMar82,WaSchuppen2000}; see \cite{YukselBasarBook} for extensions and a detailed literature review), and iterative algorithms such as sequential update laws as well as gradient based computational algorithms may converge to an optimal solution. In the context of linear quadratic Gaussian (LQG) static teams, this leads to the optimality of linear policies. 
However, characterization of convexity is in general a difficult problem and only restrictive conditions appear in the literature where these conditions do not utilize the probabilistic and information structure related aspects of convexity even for static teams. In view of this discussion, there are two main contributions of the paper:
\begin{itemize}
\item The first set of contributions involve the sets of strategic measures: Strategic measures for static and dynamic team problems will be defined and studied, and their integral representations, convexity and Borel measurability properties under the weak convergence topology will be investigated. It will be shown that static and general dynamic team problems do not lead to a convex set of strategic measures, even in the presence of common or independent randomness, unlike a small class of dynamic team problems with classical information structures (that includes fully observed single-decision-maker stochastic control problems). It will be shown that externally provided independent common randomness or private randomness is useless for static teams or dynamic team problems, respectively. Conditions for compactness of the sets of strategic measures under the weak convergence topology are presented, which lead to the existence of optimal team policies. Finally, Borel measurability properties of strategic measures and universal measurability properties of value functions are established.
\item The paper addresses the problem of when a dynamic sequential team problem is convex. It provides necessary and sufficient conditions for the convexity of static and dynamic team problems by utilizing the information structure and probability related properties. Building on these results, generalizations of Radner's \cite{rad62} and Krainak et. al.'s \cite{KraMar82}  theorems will be presented. For dynamic teams, static reduction is provided as a useful tool to establish not only convexity, but also obtain a precise method through which a probabilistic characterization of the {\it lack of a signaling incentive among decision makers} can be characterized; this generalizes the partial nestedness conditions leading to the characterization of a large class of convex teams.  
\end{itemize}

We note that our formulation in this paper considers stochastic criteria under an expected cost minimization objective. For further criteria where convexity properties have been studied, we refer the reader to \cite{Rotkowitz}, \cite{Voulgaris} and \cite{Bamieh} where operator theoretic criteria have been considered; we also refer the reader to a tutorial paper \cite{CDCTutorial} and \cite{basCDC2008}.

Here is a summary of the rest of the paper. In the following, we first provide a review of dynamic team problems, utilizing Witsenhausen's intrinsic model. Section \ref{strategic} introduces strategic measures and their representation, convexity and measurability properties, as well as implications on the structural properties and sufficient conditions for the existence of optimal solutions to team problems. Section \ref{convexStat} investigates convexity for static teams and subsection \ref{convexDynam} studies the convexity properties for dynamic teams. 

\subsection{Sequential dynamic teams and Witsenhausen's characterization of information structures}\label{witsenInfoStructureReview}

In this section, we introduce the characterizations as laid out by Witsenhausen, termed as {\it the Intrinsic Model} \cite{wit75}; see \cite{YukselBasarBook} for a more comprehensive overview and further characterizations and classifications of information structures. In this model (described in discrete time), any action applied at any given time is regarded as applied by an individual decision maker/agent, who acts only once. One advantage of this model, in addition to its generality, is that the definitions regarding information structures can be compactly described.

Suppose that in the decentralized system considered below, there is a pre-defined order in which the decision makers act. Such systems are called {\it sequential systems} (for non-sequential teams, we refer the reader to Andersland and Teneketzis \cite{AnderslandTeneketzisI}, \cite{AnderslandTeneketzisII} and Teneketzis \cite{Teneketzis2}, in addition to Witsenhausen \cite{WitsenhausenSIAM71}). Suppose that in the following, the action and measurement spaces are standard Borel spaces, that is, Borel subsets of Polish (complete, separable and metric) spaces. In the context of a sequential system, the {\it Intrinsic Model} has the following components:

\begin{itemize}
\item A collection of {\it measurable spaces} $\{(\Omega, {\cal F}),
(\mathbb{U}^i,{\cal U}^i), (\mathbb{Y}^i,{\cal Y}^i), i \in {\cal N}\}$, specifying the system's distinguishable events, and the control and measurement spaces. Here $N=|{\cal N}|$ is the number of control actions taken, and each of these actions is taken by an individual (different) DM (hence, even a DM with perfect recall can be
regarded as a separate decision maker every time it acts). The pair $(\Omega, {\cal F})$ is a
measurable space (on which an underlying probability may be defined). The pair $(\mathbb{U}^i, {\cal U}^i)$
denotes the measurable space from which the action, $u^i$, of decision maker $i$ is selected. The pair $(\mathbb{Y}^i,{\cal Y}^i)$ denotes the measurable observation/measurement space.

\item A {\it measurement constraint} which establishes the connection between the observation variables and the system's distinguishable events. The $\mathbb{Y}^i$-valued observation variables are given by $y^i=\eta^i(\omega,{\bf u}^{[1,i-1]})$, ${\bf u}^{[1,i-1]}=\{u^k, k \leq i-1\}$, $\eta^i$ measurable functions and $u^k$ denotes the action of DM $k$. Hence, the information variable $y^i$ induces a $\sigma$-field, $\sigma({\cal I}^i)$ over $\Omega \times \prod_{k=1}^{i-1} \mathbb{U}^k$
\item A {\it design constraint} which restricts the set of admissible $N$-tuple control laws $\underline{\gamma}= \{\gamma^1, \gamma^2, \dots, \gamma^N\}$, also called
{\it designs} or {\it policies}, to the set of all
measurable control functions, so that $u^i = \gamma^i(y^i)$, with $y^i=\eta^i(\omega,{\bf u}^{[1,i-1]})$, and $\gamma^i,\eta^i$ measurable functions. Let $\Gamma^i$ denote the set of all admissible policies for DM $i$ and let ${\bf \Gamma} = \prod_{k} \Gamma^k$.
\end{itemize}

We note that, the intrinsic model of Witsenhausen gives a set-theoretic characterization of information fields, however, for standard Borel spaces, the model above is equivalent to that of Witsenhausen's.

One can also introduce a fourth component:
\begin{itemize}
\item A {\it probability measure} $P$ defined on $(\Omega, {\cal F})$ which describes the measures on the random events in the model. 
\end{itemize}

Under this intrinsic model, a sequential team problem is {\it dynamic} if the
information available to at least one DM is affected by the action of at least one other DM. A decentralized problem is {\it static}, if the information available at every decision maker is only affected by exogenous disturbances (Nature); that is no other decision maker can affect the information at any given decision maker.

Information structures can also be classified as {\it classical}, {\it quasi-classical} or {\it nonclassical}. An Information Structure (IS) $\{y^i, 1 \leq i \leq N \}$ is {\it classical} if $y^i$ contains all of the information available to DM $k$ for $k < i$. An IS is {\it quasi-classical} or {\it partially nested}, if whenever $u^k$, for some $k < i$, affects $y^i$ through the measurement function $\eta^i$, $y^i$ contains $y^k$ (that is $\sigma(y^k) \subset \sigma(y^i)$). An IS which is not partially nested is {\it nonclassical}.

Let \[\underline{\gamma} = \{\gamma^1, \cdots, \gamma^N\}\]
and let a cost function be defined as:
\begin{eqnarray}\label{lossF}
J(\underline{\gamma}) = E[c(\omega_0,{\bf u})] = E[c(\omega_0,\gamma^1(y^1),\cdots,\gamma^N(y^N))],
\end{eqnarray}
for some non-negative measurable loss (or cost) function $c: \Omega_0 \times \prod_k \mathbb{U}^k \to \mathbb{R}_+$. Here, we have the notation ${\bf u}=\{u^t, t \in {\cal N}\}$. Here, $\omega_0$ may be viewed as the cost function relevant exogenous variable and is contained in $\omega$. 

\begin{defn}\label{Def:TB1}\index{Optimal team cost}
For a given stochastic team problem with a given information
structure, $\{J; \Gamma^i, i\in {\cal N}\}$, a policy (strategy) $N$-tuple
${\underline \gamma}^*:=({\gamma^1}^*,\ldots, {\gamma^N}^*)\in {\bf \Gamma}$ is
an {\it optimal team decision rule} ({\it team-optimal
decision rule} or simply {\it team-optimal solution}) if
\begin{equation}J({\underline \gamma}^*)=\inf_{{{\underline \gamma}}\in {{\bf \Gamma}}}
J({{\underline \gamma}})=:J^* , \label{eq:5}\end{equation} provided that such a strategy
exists.  The cost level achieved by this strategy, $J^*$, is the {\it
minimum {\rm (or} optimal{\rm )} team cost}.
\end{defn}
%

\begin{defn}\label{Def:TB2} \index{Person-by-person-optimality}
For a given $N$-person stochastic team with a fixed information structure, $\{J;
\Gamma^i, i \in {\cal N}\}$, an $N$-tuple of strategies
${\underline \gamma}^*:=({\gamma^1}^*,\ldots, {\gamma^N}^*)$ constitutes a {\it Nash
equilibrium} (synonymously,  a {\it person-by-person optimal} (pbp optimal) solution) if, for all $\beta
\in \Gamma^i$ and all $i\in {\cal N}$, the following inequalities hold:
\begin{equation}{J}^*:=J({\underline \gamma}^*) \leq J({\underline \gamma}^{-i*},
\beta), \label{eq:7}\end{equation} where we have adopted the notation
\begin{equation}({\underline \gamma}^{-i*},\beta):= (\gamma^{1*},\ldots, \gamma^{(i-1)*},
\beta, \gamma^{(i+1)*},\ldots, \gamma^{N*}). \label{eq:8}\end{equation}
\end{defn}

For notational simplicity, let for any $1 \leq k \leq N$, $\gamma^{-k} := \{\gamma^i, i \in \{1,\cdots,N\} \setminus \{k\} \}$

In the following, we will denote by bold letters the ensemble of random variables across the DMs; that is ${\bf y}=\{y^i, i=1,\cdots,N\}$ and ${\bf u}=\{u^i, i=1,\cdots,N\}$.

%
%
%


\subsection{Static reduction of sequential dynamic teams}\label{EquivIS}

Following Witsenhausen \cite{wit88}, we say that two information structures are equivalent if: (i) The policy spaces are equivalent/isomorphic in the loose sense that policies under one information structure are realizable under the other information structure, (ii) the costs achieved under equivalent policies are identical almost surely, and (iii) if there are constraints in the admissible policies, the isomorphism among the policy spaces preserves the constraint conditions. 

A large class of sequential team problems admit an equivalent information structure which is static. This is called the {\it static reduction} of an information structure. 

\subsubsection{Partially nested case}

An important information structure which is not nonclassical, is of the {\it quasi-classical} type, also known as {\it partially nested}\index{Partially nested information structure}; an {\sl IS} is  partially nested if an agent's information at a particular stage $t$ can depend on the action of some other agent at some stage $t' \leq t$ only if she also has access to the information of that agent at stage $t'$. For such team problems with partially nested information, one talks about {\it precedence relationships} among agents: an agent DM $i$ is {\it precedent} to another agent DM $j$ (or DM $i$  {\it communicates} to  DM $j$), if the former agent's actions affect the  information of the latter, in which case (to be partially nested) DM $j$ has to have the information based on which the action-generating policy of DM $i$ was constructed. 

 For partially nested (or quasi-classical) information structures, static reduction has been studied by Ho and Chu in the specific context of LQG systems \cite{HoChu} and for a class of non-linear systems satisfying restrictive invertibility properties \cite{ho1973equivalence}. We will discuss partially nested team problems further in Section \ref{QuasiRed}.

\subsubsection{Nonclassical case: Witsenhausen's equivalent model and static reduction of sequential dynamic teams}\label{EquivIS2}


Witsenhausen shows that a large class of sequential team problems admit an equivalent information structure which is static. This is called the {\it static reduction} of an information structure. 

An equivalence between sequential dynamics teams and their static reduction is as follows (termed as {\it the equivalent model} \cite{wit88}\index{Witsenhausen's Equivalent Model}).

Consider a dynamic team setting according to the intrinsic model where there are $N$ time stages, and each DM observes, for some $t=1,2,\cdots,N$, \[y^t=\eta_t(\omega,u^1,u^2,\cdots,u^{t-1}),\] and the decisions are generated by $u^t=\gamma_t(y^t)$. 
Here, as before, $\omega$ is the collection of primitive (exogenous) variables. 
The resulting cost under a given team policy is, as in (\ref{lossF})
\[J(\underline{\gamma}) = E[c(\omega_0,{\bf u})].\]
This dynamic team can be converted to a static team provided that the following absolute continuity condition holds: For every $t \in {\cal N}$, there exists a function $f_t$ such that for all Borel $S \subset \mathbb{Y}^t$:
\begin{eqnarray}\label{staticReduc}
&& P(y^t \in S | \omega,u^1,\cdots,u^{t-1}) \nonumber \\
&& \quad \quad \quad \quad = \int_{S} f_t(\omega,u^1,u^2,\cdots,u^{t-1},y^t) Q_t(dy^t), 
\end{eqnarray}
where $Q_t$ is some arbitrary reference probability measure for the variable $y^t$. 
Under any fixed team policy, we can then write
\[P(d\omega,d{\bf y}) = P(d\omega) \prod_{t=1}^N f_t(\omega,u^1,u^2,\cdots,u^{t-1},y^t) Q_t(dy^t).\]
The cost function $J(\underline{\gamma})$ can then be written as
\begin{eqnarray}\label{staticReduc2}
J(\underline{\gamma})= \int P(d\omega) \prod_{t=1}^N (f_t(y_t,\omega,u^1,u^2,\cdots,u^{t-1},y^t) Q_t(dy^t))c(\omega_0,{\bf u}),
\end{eqnarray}
where now the measurement variables can be regarded as independent and by incorporating the $\{f_t\}$ terms into $c$, we can obtain an equivalent {\it static team} problem. Hence, the essential step is to appropriately adjust the probability space and the cost function. The new cost function may now explicitly depend on the measurement values, such that
\[c_s(\omega,{\bf y}, {\bf u}) = c(\omega_0,{\bf u}) \prod_{t=1}^N f_t(y_t,\omega,u^1,u^2,\cdots,u^{t-1},y^t).\]
In this case, we can view $\omega, {\bf y}$ as the cost-relevant exogenous variable: By an abuse of notation, we will use the same notation $\omega_0$ to denote $\omega, {\bf y}$ when it is clear that such a cost function comes from a static reduction.

We note that, as Witsenhausen notes in \cite{wit88}, a static reduction always holds when the measurement variables take values from countable set since a reference measure as in $Q_t$ above can be constructed on the measurement variable $y^t$ (e.g., $Q_t(y^t) = \sum_{i \geq 1} 2^{-i} 1_{\{y^t = m_i\}}$ where $\mathbb{Y}^t=\{m_i, i \in \mathbb{N}\}$) so that the absolute continuity condition always holds. On the other hand, for continuous spaces, observe that under a control-sharing pattern with $y^2=u^1$, the absolute continuity condition required for Witsenhausen's static reduction may fail: $P(y^2 \in A|u^1)=1_{\{u^1 \in A\}}$, leading to a delta function supported at $u^1$ and if the reference measure $\mu$ with $y^2 \sim \mu$ admits a density, the absolute continuity condition will not hold. We also note that a continuous-time generalization for static reduction similar to Girsanov's method has been presented by Charalambous and Ahmed \cite{charalambous2014equivalence}.

\subsection{Convex static teams: Radner's and Krainak et.al.'s theorems}

An important property of convexity is that local optimality conditions imply global optimality conditions, as we briefly detail below. We note that a more general characterization of convexity will be presented later in the paper. 

\begin{defn}\label{Def:TB4}
Given a static stochastic team problem $\{J; \Gamma^i,i\in\IN\}$, a
policy $N$-tuple $\underline{\gamma}\in \boldGamma$ is {\it
stationary} if (i) $J(\underline{\gamma})$ is finite, (ii)~the $N$
partial derivatives in the following equations are well defined, and (iii)~$\underline{\gamma}$ satisfies these equations:
\begin{equation}\left[ \nabla_{u^i} {E_{\omega | y^i}}
c(\omega_0;{\underline{\gamma}}^{-i}(\boldy), u^i) \right]
|_{u^i=\gamma^i(y^i)} = 0, a.s.\quad i\in\IN .\label{eq:19}\end{equation} \index{Stationarity condition}
\end{defn}



The following results are due to Krainak et. al. \cite{KraMar82} and \cite{YukselBasarBook}, generalizing Radner \cite{rad62}. We follow the presentation in \cite{YukselBasarBook}, which also contains the proofs of the results.

\begin{thm}\label{Thm:TB4}\cite{rad62} \cite{KraMar82}
Let $\{ J; \Gamma^i , i\in\IN \}$ be a static stochastic team
problem where $\mathbb{U}^i = \mathbb{R}^{m_i}, i\in\IN$, the loss function $c
( \omega_0 , \boldu )$ is convex and continuously differentiable
in $\boldu$ a.s., and $J (\undergamma)$ is bounded from below on $\boldGamma$.
Let $\undergamma^* \in \boldGamma$ be a policy $N$-tuple with a
finite cost $(J(\undergamma^*) < \infty)$, and suppose that for
every $\undergamma \in \boldGamma$ such that $J(\undergamma) <
\infty$,
 the following holds: 
 \begin{equation} 
 \sum_{i \in {\cal N}}E \{ \nabla_{u^i} c(\omega_0 ; \undergamma^* (\boldy)) [ \gamma^i (y^i)
- \gamma^{i*} (y^i)] \} \geq 0, 
\label{eq:200}
\end{equation} 
where $E \{ \cdot \}$ denotes the total expectation.  
Then,
$\undergamma^*$ is a {\it team-optimal policy}, and it is {\it
unique} if $c$ is strictly convex in $\boldu$.
\end{thm}

Note that the conditions of {\sl Theorem~\ref{Thm:TB4}} above do not include the
stationarity of $\undergamma^*$, and furthermore inequalities~(\ref{eq:200})
may not generally be easy to check, since they involve all permissible
policies $\undergamma$ (with finite cost). Instead, either one of the following two conditions will achieve this objective \cite{KraMar82} \cite{YukselBasarBook}:
\begin{description}
\item{(c.5)} For all $\undergamma \in \Gamma$ such that $J
(\undergamma) < \infty$, the following random variables have
well-defined (finite) expectations
$$
\nabla_{u^i} c(\omega_0; \undergamma^* (\boldy)) [ \gamma^i (y^i) -
\gamma^{i*} (y^i) ], \quad i \in {\cal N}
$$
\item{(c.6)} $\Gamma^i$ is a Hilbert space for each $i\in\IN$, and
$J(\undergamma) < \infty$ for all $\undergamma \in \Gamma$.
Furthermore, $$ E_{\omega | y^i} \{ \nabla_{u^i} c(\omega_0 ; \undergamma^*
(\boldy) \} \in
\Gamma^i, \quad i\in\IN.
$$
\end{description}

\begin{thm}\label{Thm:TB5}\cite{KraMar82} \cite{YukselBasarBook}
Let $\{ J; \Gamma^i, i\in\IN \}$ be a static
stochastic team problem which satisfies all the hypotheses of
Theorem~\ref{Thm:TB4}, with the exception of the inequality (\ref{eq:200}).  Instead
of~(\ref{eq:200}), let either (c.5) or (c.6) be satisfied.  Then, if
$\undergamma^* \in \boldGamma$ is a stationary policy it is also
team optimal.  Such a policy is unique if $c(\omega_0 ; \boldu )$ is
strictly convex in $\boldu$, a.s.
\end{thm}

\section{Strategic Measures, Convexity Properties, and Optimal Solutions}\label{strategic}

\subsection{Viewing measurable policies as a subset of randomized policies and strategic measures}

We can view a measurable policy as a special case of randomized policies. This interpretation has many useful properties, one being the topological use of the space of probability measures under weak convergence. We recall here the following representation result due to Borkar \cite{BorkarRealization}. Let $\mathbb{X}, \mathbb{M}$ be standard Borel spaces. Let the notation ${\cal P}(\mathbb{X})$ denote the set of probability measures on $\mathbb{X}$. Consider the set of probability measures
\begin{eqnarray}\label{extremePointQuan0}
\Theta: = \{\zeta \in {\cal P}(\mathbb{X} \times \mathbb{M}): \zeta(dx,dm) = P(dx)Q^f(dm|x), Q^f(\cdot | x) = 1_{\{f(x) \in \cdot\}}, f : \mathbb{X} \to \mathbb{M} \}, \nonumber 
\end{eqnarray}
on $\mathbb{X} \times \mathbb{M}$ having fixed input marginal $P$, equipped with weak topology. This set is the (Borel measurable) set of the extreme points of the set of probability measures on $\mathbb{X} \times \mathbb{M}$ with a fixed input marginal $P$. For compact $\mathbb{M}$, the Borel measurability of $\Theta$ follows from \cite{Choquet} since set of probability measures on $\mathbb{X} \times \mathbb{M}$ with a fixed input marginal $P$ is a convex and compact set in a complete separable metric space, and therefore, the set of its extreme points is Borel measurable.  But the non-compact case also holds; see Lemma 2.3 in \cite{BorkarRealization}. Furthermore, given a fixed input marginal measure $P$ on $\mathbb{X}$, any stochastic kernel from $\mathbb{X}$ to $\mathbb{M}$ can be identified by a probability measure ${\cal K} \in {\cal P}(\mathbb{X} \times \mathbb{M})$ such that
\begin{eqnarray}\label{convR0}
{\cal K} (A)  = \int_{\Theta} \xi(dQ) Q(A), \quad  A \in {\cal B}(\mathbb{X} \times \mathbb{M})
\end{eqnarray}
for some $\xi \in {\cal P}(\Theta)$. In particular, a stochastic kernel can thus be viewed as an integral representation over probability measures induced by deterministic policies. 

%

For a team setup, for any DM $k$, let
\begin{eqnarray}\label{extremePointQuan1}
\Theta^k: &=& \bigg\{\zeta \in {\cal P}(\mathbb{Y}^k \times \mathbb{U}^k): \zeta = P_kQ^{\gamma^k}, Q^{\gamma^k}(\cdot|y^k) = 1_{\{\gamma^k(y^k) \in \cdot\}}, \gamma^k \in \Gamma^k, \nonumber \\
&& \quad \quad \quad \quad \quad \quad \quad \quad \quad \quad \quad \quad \quad\quad \quad P_k(\cdot) = P(y^k \in \cdot) \bigg\} \nonumber 
\end{eqnarray}
For a static team, $P_k$ would be independent of the policies of the preceding DMs. The set of all stochastic kernels from $\mathbb{Y}^k$ to $\mathbb{U}^k$ with fixed input marginal measure $P_k$ is such that any element ${\cal K}^k$ in this space can be expressed in the  form
\begin{eqnarray}\label{convR}
{\cal K}^k (A)  = \int_{\Theta^k} \xi^k(dQ) Q(A), \quad  A \in {\cal B}(\mathbb{Y}^k \times \mathbb{U}^k)
\end{eqnarray}
for some $\xi \in {\cal P}(\Theta^k)$. 

For stochastic control problems, {\it strategic measures} are defined (see \cite{Schal}, \cite{piunovskii1998controlled}, \cite{dynkin1979controlled} and \cite{feinberg1996measurability}) as the set of probability measures induced on the product spaces of the state and action pairs by measurable control policies: Given an initial distribution on the state, and a policy, one can uniquely define a probability measure on the product space. Certain measurability, compactness and convexity properties of strategic measures for single decision maker problems were studied in \cite{dynkin1979controlled}, \cite{piunovskii1998controlled}, \cite{feinberg1996measurability} and \cite{blackwell1976stochastic}.  In the following, we discuss the case for team problems, study some convexity properties and implications on the existence and structure of optimal policies.

\subsection{Convexity properties of sets of strategic measures and redundancy of common or private independent randomness}

Consider a static team problem defined under Witsenhausen's intrinsic model in Section \ref{witsenInfoStructureReview}. Let $L_A(\mu)$ be the set of strategic measures induced by all admissible team policies with $(\omega_0, {\bf y}) \sim \mu$. In the following, $B = B^0 \times \prod_{k} (A^k \times B^k)$ are used to denote the Borel sets in $\Omega_0 \times \prod_k (\mathbb{Y}^k \times \mathbb{U}^k)$,
\begin{eqnarray}
L_A(\mu)&:=& \bigg\{P \in {\cal P}\bigg(\Omega_0 \times \prod_{k=1}^N (\mathbb{Y}^k \times \mathbb{U}^k)\bigg): \nonumber \\
&& \quad P(B) = \int_{B^0 \times \prod_k A^k} \mu(d\omega_0, d{\bf y}) \prod_k 1_{\{u^k = \gamma^k(y^k) \in B^k\}},  \nonumber \\
&& \quad \quad \quad \quad \gamma^k \in \Gamma^k, B \in {\cal B}(\Omega_0 \times \prod_k (\mathbb{Y}^k \times \mathbb{U}^k)) \bigg\}
\end{eqnarray}
Let $L_A(\mu,\underline{\gamma})$ be the strategic measure under a particular $\underline{\gamma} \in {\bf \Gamma}$. 
Let $L_R(\mu)$ be the set of strategic measures induced by all admissible team policies where $\omega, {\bf y} \sim \mu$ and policies are individually randomized (that is, with independent randomizations):
\[L_R(\mu) := \bigg\{P \in {\cal P}\bigg(\Omega_0 \times \prod_{k=1}^N (\mathbb{Y}^k \times \mathbb{U}^k)\bigg): P(B) = \int_{B} \mu(d\omega_0, d{\bf y}) \prod_k \Pi^k(d u^k| y^k)\bigg\}\]
where $\Pi^k$ takes place from the set of stochastic kernels from $\mathbb{Y}^k$ to $\mathbb{U}^k$ for each $k$.
Consider $\Upsilon = [0,1]^{N}$ and let 
\begin{eqnarray}\label{LC}
L_{C}(\mu) := \bigg\{P \in {\cal P}\bigg(\Omega_0 \times \prod_{k=1}^N (\mathbb{Y}^k \times \mathbb{U}^k)\bigg): P(B) = \int \eta(dz) L_A(\mu,\underline{\gamma}(z))(B), \quad \eta \in{\cal P}(\Upsilon) \bigg\}  \nonumber \\
\end{eqnarray}
where $\underline{\gamma}(z)$ denotes a collection of team policies measurably parametrized by $z \in \Upsilon$ so that the map $L_A(\mu,\underline{\gamma}(z)): \Upsilon \to L_A(\mu)$ is Borel measurable.
Finally, let $L_{CR}$ denote the set of strategic measures that are induced by some fixed but common independent randomness and arbitrary private independent randomness.
\begin{eqnarray}
L_{CR}(\mu)&:=& \bigg\{P \in {\cal P}\bigg(\Omega_0 \times \prod_{k=1}^N (\mathbb{Y}^k \times \mathbb{U}^k)\bigg): \nonumber \\
&&  \quad \quad \quad \quad P(B) = \int_{B \times [0,1]^N} \eta(dz) \mu(d\omega_0, d{\bf y}) \prod_k \Pi^k(d u^k| y^k,z) \bigg\} \nonumber 
\end{eqnarray}
\begin{thm}\label{FeinbergStrategic}
Consider a static team.
\begin{itemize}
\item (i) $L_R$ has the following representation
\begin{eqnarray}\label{representLR}
L_{R}(\mu)= &&\{P \in {\cal P}\bigg(\Omega_0 \times \prod_{k=1}^N (\mathbb{Y}^k \times \mathbb{U}^k)\bigg): P(B) = \int U(dz) L_A(\mu,\underline{\gamma}(z))(B), \nonumber \\
&&   U \in{\cal P}(\Upsilon), U(dv_1,\cdots,dv_N) = \prod_s \eta_k(dv_k), \eta_k \in {\cal P}([0,1]) \},
\end{eqnarray}
that is $U \in {\cal P}(\Upsilon)$ is constructed by the product of $N$ independent random variables on $[0,1]$.
\item (ii)  $L_{C}(\mu)$ is convex. Its extreme points form $L_A(\mu)$. Furthermore, $L_R(\mu) \subset L_C(\mu)$.
\item (iii) \[\inf_{\underline{\gamma} \in {\bf \Gamma}} J(\underline{\gamma}) = \inf_{P \in L_A(\mu)} \int P(ds) c(s) = \inf_{P \in L_R(\mu)} \int P(ds) c(s) = \inf_{P \in L_{C}(\mu)} \int P(ds) c(s)\]
In particular, deterministic policies are optimal among the randomized class.
\item (iv) The sets $L_R(\mu)$ and $L_{CR}(\mu)$ are not convex. In particular, the presence of independent or common randomness does not convexify the set of strategic measures.
\end{itemize}
\end{thm}


\textbf{Proof.}
\begin{itemize}
\item (i) Any stochastic kernel from a Polish space to another one, $P(dx|y)$, can be realized by some measurable function $x=f(y,v)$ where $v$ is a $[0,1]$-valued independent random variable and $f$ is measurable (see Lemma 1.2 in Gikhman and Shorodhod \cite{gihman2012controlled}, Theorem 1 in Feinberg \cite{Feinberg1} or Lemma 3.1 of Borkar \cite{BorkarRealization}). In particular, in the representation (\ref{convR}), since $\Theta^k$ is a Borel set in a Polish space, by the Borel isomorphism theorem (Appendix 1 in \cite{dynkin1979controlled} or  {\sl Chapter 13} in \cite{Dud02}), there exists a bijection $\kappa$ between $[0,1]$ and $\Theta^k$ so that ${\cal K}^k (A)  = \int p(ds) P^{\kappa(s)}(A)$, for some $p \in {\cal P}([0,1])$ defined by $p(\kappa^{-1}(A))=\xi(A)$. 
The result then follows by replacing each of the stochastic kernels with such a representation due to \cite{gihman2012controlled}. 


\item (ii) For convexity, observe that for $\kappa_1 \in L_C(\mu)$, $\kappa_2 \in L_C(\mu)$, and $\alpha \in (0,1)$, $\alpha \kappa_1 + (1 - \alpha) \kappa_2 \in L_C(\mu)$ since with
\[ \kappa_1(B) = \int \eta_1(dz) L_A(\mu,\underline{\gamma}(z)), \quad \eta \in{\cal P}(\Upsilon) \]
\[ \kappa_2(B) = \int \eta_2(dz) L_A(\mu,\underline{\gamma}(z)), \quad \eta \in{\cal P}(\Upsilon), \]
\[ (\alpha \kappa_1 + (1 - \alpha) \kappa_2)(B) = \int (\alpha \eta_1 + (1 - \alpha) \eta_2)(dz) L_A(\mu,\underline{\gamma}(z)), \quad \eta \in{\cal P}(\Upsilon) \]

For the extreme sets being in $L_A(\mu)$: Suppose $\kappa =  \int \eta(dz) L_A(\mu,\underline{\gamma}(z))$ is not in $L_A(\mu)$, then it can't be extreme since for every such $\kappa$ one can construct $\kappa_1$, $\kappa_2$ and $\lambda \in (0,1)$ so that $\kappa = \lambda \kappa_1 + (1 - \lambda) \kappa_2$: Partition the support set of $\kappa$ to two disjoint components $D_1, D_2 \in {\cal B}(\Upsilon)$ with $\lambda = \eta(D_1)$ and define $\kappa_1 = \int_{D_2} {\eta(dz) \over \lambda} L_A(\mu,\underline{\gamma}(z))$, $\kappa_2 = \int_{D_1} {\eta(dz) \over 1- \lambda} L_A(\mu,\underline{\gamma}(z))$. For the converse, suppose $\kappa$ is not extreme, then there exists $\alpha, \eta_1, \eta_2$ such that $\kappa = \int (\alpha \eta_1 + (1 - \alpha) \eta_2)(dz) L_A(\mu,\underline{\gamma}(z))$ and in particular for some set of $z$ values with a positive measure, the measure cannot be represented as: $\delta_{z} \prod_k 1_{\{\gamma^k(y^k,z) \in \cdot\}}$ for some Dirac $\delta_{z}$ measure: If this were possible, $\eta = \alpha \eta_1 + (1 - \alpha) \eta_2$,  since this has to be a mixture of two such probability measures. 

That, $L_R(\mu) \subset L_C(\mu)$ follows from the representation in (\ref{representLR}) and the definition in (\ref{LC}).

\item (iii) The inequalities
\[\inf_{\underline{\gamma} \in {\bf \Gamma}} J(\underline{\gamma}) = \inf_{P \in L_A(\mu)} \int P(ds) c(s) \geq \inf_{P \in L_R(\mu)} \int P(ds) c(s) \geq \inf_{P \in L_{C}(\mu)} \int P(ds) c(s)\]
follow since $L_A(\mu) \subset L_R(\mu) \subset L_C(\mu)$. On the other hand, the integral $\int P(ds) c(s)$ is a linear mapping from the convex set $L_C(\mu)$ to $\mathbb{R}$ (possibly to be extended to include $+ \infty$). The extreme points of $L_{C}(\mu)$ are those that are in $L_A(\mu)$. It then follows from convex analysis (see e.g. Section 2.3 in \cite{piunovskii1998controlled}) that
\[\inf_{P \in L_{C}(\mu)} \int P(ds) c(s) = \inf_{\underline{\gamma} \in {\bf \Gamma}} J(\underline{\gamma}).\]
\item (iv) A convex combination of two measures leads to a new measure 
 which will not necessarily be the product measure of two marginal probability measures: Consider two sets of measures $P^1 P^2$ and $\bar{P}^1 {\bar P}^2$. For $\lambda \in (0,1)$, a convex combination $\lambda P^1 P^2 + (1 - \lambda) \bar{P}^1 {\bar P}^2$ cannot be written as the product of their marginals $(\lambda P^1 + (1 - \lambda) \bar{P}^1) ( \lambda P^2  + {1 - \lambda} {\bar P}^2)$. 

\end{itemize}

\qed

We state two remarks in the following. 

\begin{remark}
Consider the following set of probability measures: 
\[L_{M}(\mu)= \{P: P(d\omega_0, {\bf y}) = \mu(d\omega_0, d{\bf y}), \quad P(du^i | y^i, {\bf y}^{-i}) = P(du^i|y^i), \quad \quad i \in {\cal N}\}\]
This set is a convex set and arises in applications in information theory in the context of converse theorems in multi-terminal source coding (see e.g. \cite{wagner2008improved}). It is tempting to claim that $L_M(\mu) = L_C(\mu)$, however, a counterexample in \cite{anantharam2007common} reveals that this is not the case. 
\end{remark}

\begin{remark}In the representation of $L_R(\mu)$ above in Theorem \ref{FeinbergStrategic}-(i), we could also have gone one step further and allowed for the probability measure on $[0,1]$ to be the Lebesgue (uniform) measure; but this would require us to alter the set $L_A(\mu)$ in the way the randomization appears. We will revisit this in the following (where we build on a relevant result due to Gikhman and Skorodhod \cite{gihman2012controlled}), where we consider sequential dynamic teams, instead of only static teams. 
\end{remark}


We present the following characterization for strategic measures in dynamic sequential teams. Let for all $n \in \mathbb{N}$, $h_n = \{\omega_0,y^1,u^1,\cdots,y^{n-1},u^{n-1},y^n,u^n\}$, and $p_n(dy^n|h_{n-1}) = p(dy^n | \omega_0,u^{1},\cdots,u^{n-1})$ be the transition kernel characterizing the measurements of DM $n$ according to the intrinsic model. We note that this may be obtained by the relation:
\begin{eqnarray}
p_n(y^n \in \cdot | \omega_0,u^{1},\cdots,u^{n-1}) := P\bigg(\eta^n(\omega, u^{1},\cdots,u^{n-1}) \in \cdot  \bigg| \omega_0,u^{1},\cdots,u^{n-1}\bigg) \label{kernelDefn}
\end{eqnarray}

Let $L_A(\mu)$ be the set of strategic measures induced by deterministic policies and let $L_R(\mu)$ be the set of strategic measures induced by independently provided randomized policies. We note as earlier that such an individual randomized policy can be represented in a functional form: By Lemma 1.2 in Gikhman and Shorodhod \cite{gihman2012controlled} and Theorem 1 in Feinberg \cite{Feinberg1}, for any stochastic kernel $\Pi^k$ from $\mathbb{Y}^k$ to $\mathbb{U}^k$, there exists a measurable function $\gamma^k: [0,1] \times \mathbb{Y}^k \to \mathbb{U}^k$ such that
\begin{eqnarray}\label{GSLemma}
m\{r: \gamma^k(r,y^k) \in A\} = \Pi^k(u^k \in A | y^k),
\end{eqnarray}
and $m$ is the uniform distribution (Lebesgue measure) on $[0,1]$. 

\begin{thm}\label{StrategicCharacterization}
\begin{itemize}
\item A probability measure $P \in {\cal P}\bigg(\Omega_0 \times \prod_{k=1}^N (\mathbb{Y}^k \times \mathbb{U}^k)\bigg)$ is a strategic measure induced by a deterministic policy (that is in $L_A(\mu)$) if and only if for every $n \in \{1,\cdots,N\}$:
\[\int P(dh_{n-1},dy^n) g(h_{n-1},y^{n}) = \int P(dh_{n-1}) \bigg(\int_{\mathbb{Y}^n} g(h_{n-1},z) p_n(dz | h_{n-1}) \bigg),\]
and
\[\int P(dh_n) g(h_{n-1},y^n,u^{n}) = \int P(dh_{n-1},dy^n) \bigg(\int_{\mathbb{U}^n} g(h_{n-1},y^{n},a) 1_{\{\gamma^n(y^n) \in da\}} \bigg),\]
for some $\gamma^n \in \Gamma^n$, for all continuous and bounded $g$, with $P(d\omega_0) = \mu(dw_0)$. 
\item A probability measure $P \in {\cal P}\bigg(\Omega_0 \times \prod_{k=1}^N (\mathbb{Y}^k \times \mathbb{U}^k)\bigg)$  is a strategic measure induced by a randomized policy (that is in $L_R(\mu)$)  if and only if for every $n \in \{1,\cdots,N\}$:
\begin{eqnarray}\label{convD1}
\int P(dh_{n-1},dy^n) g(h_{n-1},y^{n}) = \int P(dh_{n-1}) \bigg(\int_{\mathbb{Y}^n} g(h_{n-1},z) p_n(dz | h_{n-1}) \bigg), \nonumber \\
\end{eqnarray}
and
\begin{eqnarray}\label{convD2}
\int P(dh_n) g(h_{n-1},y^n,u^{n}) = \int P(dh_{n-1},dy^n) \bigg(\int_{\mathbb{U}^n} g(h_{n-1},y^{n},a) \Pi^n(da | y^n) \bigg), \nonumber \\
\end{eqnarray}
for some stochastic kernel $\Pi^n$ on $\mathbb{U}^n$ given $\mathbb{Y}^n$, for all continuous and bounded $g$, with $P(d\omega_0) = \mu(dw_0)$. 
\end{itemize}
\end{thm}

\textbf{Proof.} 
The proof follows from the fact that testing the equalities such as (\ref{convD1}-\ref{convD2}) on continuous and bounded functions implies this property for any measurable and bounded function (that is, continuous and bounded functions form a {\it separating class}, see e.g. p. 13 in \cite{Billingsley} or Theorem 3.4.5 in \cite{ethier2009markov}) and thus we recover the measurability properties leading to Witsenhausen's intrinsic model.
\qed

\begin{remark}
Theorem \ref{StrategicCharacterization} will be useful later when we will establish sufficient conditions for the compactness and Borel measurability properties (under the weak topology) of the set of strategic measures. This will then lead to the existence of optimal team policies in Sections \ref{exis} and \ref{furtherExis}, and further Borel and universal measurability properties for sets of strategic measures and value functions.
\end{remark}

\begin{thm}\label{FeinbergStrategic2}
Consider a sequential (static or dynamic) team. 
\begin{itemize}
\item(i) $L_R(\mu)$ has the following representation so that for any $P \in L_{R}(\mu)$, with $B = B^0 \times \prod (A^k \times B^k) \in {\cal B}(\Omega_0 \times \prod_{k=1}^N (\mathbb{Y}^k \times \mathbb{U}^k))$,
\begin{eqnarray}\label{representLR2}
P(B) &=& \int U(dz) L_A(\mu,\underline{\gamma}(z))(B) \nonumber \\
&& \quad \quad \quad U(dv_1,\cdots,dv_N) = \prod_s \eta (dv_s), \eta_s \in {\cal P}([0,1]),
\end{eqnarray}
where $\eta(a,b)=b-a, 0 \leq a \leq b \leq 1$, that is $\eta$ is the Lebesgue measure on $[0,1]$ and $\underline{\gamma}(z)$ is a collection of deterministic policies parametrized by $z$.
\item (ii) \[\inf_{\underline{\gamma} \in {\bf \Gamma}} J(\underline{\gamma}) = \inf_{P \in L_A(\mu)} \int P(ds) c(s) = \inf_{P \in L_R(\mu)} \int P(ds) c(s)\]
In particular, deterministic policies are optimal among the randomized class.
\end{itemize}
\end{thm}
\textbf{Proof.} 
\begin{itemize} 
\item(i) Here, we build on Lemma 1.2 in Gikhman and Shorodhod \cite{gihman2012controlled} discussed earlier in (\ref{GSLemma}). For any stochastic kernel $\pi^k$ from $\mathbb{Y}^k$ to $\mathbb{U}^k$, there exists a measurable function $\gamma^k: [0,1] \times \mathbb{Y}^k \to \mathbb{U}^k$ such that $ m\{r: \gamma^k(r,y^k) \in A\} = \Pi^k(u^k \in A | y^k)$,
and $m$ is the uniform (Lebesgue measure) on $[0,1]$. By concatenating $N$ such measures in a product form on $\Upsilon=[0,1]^N$, we obtain the representation result.
\item(ii) This follows as in Theorem \ref{FeinbergStrategic} (iii), but rather than a convex analytic argument, here we use the integral representation as a consequence of (i); see also Theorem 1-c in \cite{Feinberg1}. 
\end{itemize}
\qed

If a sequential team is not classical (and not necessarily {\it non-classical}), the set of strategic measures is not convex. Note here that the team can be static or dynamic. For the special case of single-decision maker stochastic control problems with fully available state information at the controller (note that this can be viewed as a team problem with a classical information structure), convexity of the set of strategic measures has been established in \cite{dynkin1979controlled} and \cite{piunovskiy1997multicriteria}. We present a concise proof for the more general team problems with classical information structures.

\begin{thm} [Convexity for classical information structures] If the information structure is classical and the information structure is expanded so that when independently randomized policies are allowed DM $i$ has access to $y^k, u^k, k < i$ and $y^i$, then the set of strategic measures is convex. 
\end{thm}

\textbf{Proof.}
The proof follows from Theorem \ref{StrategicCharacterization} and the fact that the condition on the control policies is redundant. Nonetheless, to present a direct proof, note that with $N=1$, convexity is immediate since any stochastic kernel can be realized by a representation of the form $\gamma^1(y^1,r^1)$ for an independent $[0,1]$ valued $r^1$. With $N > 1$, consider two strategic measures defined sequentially for $1 \leq k \leq N$,
\begin{eqnarray}
&& P^1_k(d{\bf u}, d{\bf y},d\omega_0) = \Pi^1(du^{k}|y^1,\cdots,y^k,u^1,\cdots,u^{k-1}) P(dy^k | y^1,\cdots,y^{k-1},u^1,\cdots,u^{k-1},\omega_0) \nonumber \\
&& \quad \quad \quad \quad \quad \quad \quad \quad \times  P^1_{k-1}(y^1,\cdots,y^{k-1},u^1,\cdots,u^{k-1},d\omega_0)
\end{eqnarray}
\begin{eqnarray}
&& P^2_k(d{\bf u}, d{\bf y},d\omega_0) = \Pi^2(du^{k}|y^1,\cdots,y^k,u^1,\cdots,u^{k-1}) P(dy^k | y^1,\cdots,y^{k-1},u^1,\cdots,u^{k-1},\omega_0) \nonumber \\
&& \quad \quad \quad \quad \quad \quad \quad \quad \times P^2_{k-1}(y^1,\cdots,y^{k-1},u^1,\cdots,u^{k-1},d\omega_0)
\end{eqnarray}
Now, write the above as:
\[P^1(d{\bf u}, d{\bf y},d\omega_0) = P^1(y^1,\cdots,y^k,u^1,\cdots,u^{k} | u^1, y^1) P^1(du^1|y^1) \mu(dy^1,d\omega_0)\]
\[P^2(d{\bf u}, d{\bf y},d\omega_0) = P^2(y^1,\cdots,y^k,u^1,\cdots,u^{k} | u^1, y^1) P^2(du^1|y^1) \mu(dy^1,d\omega_0)\]
We would like to see if $\lambda P^1 + (1 - \lambda) P^2$ can be realized. If one takes the convex combinations $\lambda P^1(du^1|y^1) + (1 - \lambda) P^2(du^1|y^1)$, the convex combination will be realized by not altering any of the policies adopted by the subsequent DMs.  \qed

An implication of the representation results in Theorems \ref{FeinbergStrategic} and \ref{FeinbergStrategic2} is that an optimal policy can be obtained within the deterministic team policies. Earlier, we used Blackwell's Principle of Irrelevant Information (\cite{Blackwell2} \cite{Blackwell3}; see \cite[p. 457]{YukselBasarBook}) to derive such a result in a setting where the randomizations were restricted to be independent (see \cite{gupta2014existence} for further discussions and {\sl Chapter 4} of \cite{YukselBasarBook} in the context of comparison of experiments). 

\begin{thm}\label{BIR}\cite{gupta2014existence}
For a sequential dynamic team, any randomized team policy (with independent and private randomization) can be replaced with a policy which is deterministic and which performes at least as good as the original (randomized) policy. 
\end{thm}


Both results have their own merits: (i) Theorem \ref{BIR} allows one to replace any sequential team policy with independently randomized policies with a deterministic one. However, it requires {\it private independent randomness}, that is {\it common randomness} is not allowed.  On the other hand, (ii) Theorems \ref{FeinbergStrategic} and \ref{FeinbergStrategic2} use the powerful theory of convex analysis; they allow for establishing the fact that even {\it common independent randomness} does not help for static teams or {\it private randomness} does not help for sequential teams; these theorems do not utilize backwards induction, and are thus applicable to teams with countably many DMs (by generalizing the representation result accordingly) as well as allow for the analysis to be applicable to constrained team problems. 


\subsection{Existence of optimal team policies}\label{exis}

Establishing the existence and structure of optimal policies is a challenging problem. Existence of optimal policies for static and a class of sequential dynamic teams have been studied recently in \cite{gupta2014existence}. More specific setups and non-existence results have been studied in \cite{WuVer11}, \cite{wit68}, \cite{YukselOptimizationofChannels} and \cite{YukselBasarBook}. Existence of optimal team policies has been established in \cite{charalambous2014maximum} for a class of continuous-time decentralized stochastic control problems. For a class of teams which are {\it convex}, one can reduce the search space to a smaller parametric class of policies, as discussed earlier. Finally, the strategic measure approach for single-decision maker problems and fully observed Markov models has been studied in \cite{Schal} and \cite{piunovskii1998controlled}, among other contributions in the literature. 

For team problems, considering the set of strategic measures and compactification or convexification of these sets of measures through introducing private or common randomness will allow for placing a useful topology, that of weak convergence of probability measures, on the strategy spaces. Combined with a compactness condition, this allows for establishing the existence of optimal team policies.

%
%
%

We recall that a sequence of probability measures $\mu_n$ on a standard Borel space converges to some probability measure $\mu$ weakly if $\int \mu_n(dz) f(z) \to \int \mu(dz) f(z)$ for all continuous and bounded $f$. The sequence converges setwise if $\int \mu_n(dz) f(z) \to \int \mu(dz) f(z)$ for all measurable and bounded $f$. Thus, setwise convergence is stronger than weak convergence.

\begin{thm}\label{existenceT}
(i) Consider a static or dynamic team. Let the loss function $c$ be lower semi-continuous in $(\omega_0, {\bf u})$ and $L_R(\mu)$ be a compact subset under weak topology. Then, there exists an optimal team policy. This policy is deterministic and hence induces a strategic measure in $L_A(\mu)$.  \\
(ii) Consider a static team or the static reduction of a dynamic team with $c$ denoting the loss function. Let $c$ be lower semi-continuous in $\omega_0, {\bf u}$ and $L_C(\mu)$ be a compact subset under weak topology. Then, there exists an optimal team policy. This policy is deterministic and hence induces a strategic measure in $L_A(\mu)$.  \\
\end{thm}

\textbf{Proof.} The proofs of both results build on the result that for a lower semi-continuous (and not necessarily bounded) $f: \mathbb{X} \to \mathbb{R}$, where $\mathbb{X}$ is Polish; $\int f(x) \mu(dx)$ is lower semi-continuous on ${\cal P}(\mathbb{X})$ under weak convergence (see, e.g., Lemma 4.3 in \cite{villani2008optimal}). Weierstrass theorem then leads to the existence of an optimal strategic measure, and thus an optimal policy (which is deterministic as a result of Theorems \ref{FeinbergStrategic}) or \ref{BIR}.  
\qed

In the following, we present and review some sufficient conditions for the compactness of $L_R(\mu)$ under weak topology. However, we first show that unless certain conditions are imposed, the conditional independence property is not preserved under weak or setwise convergence and thus $L_R(\mu)$ is in general not compact.

\begin{thm}
Consider a sequence of probability measures $P_n \in {\cal P}(\mathbb{U}^1 \times \mathbb{Y} \times \mathbb{U}^2)$ so that for all $n$: 
\[P_n(du^1| y, u^2) = P_n(du^1|y).\]
If $P_n \to P$ setwise (and thus also weakly), it is not necessarily the case that 
\[P(du^1| y, u^2) = P(du^1|y).\]
That is, conditional independence of $u^1$ and $u^2$ given $y$ is not preserved under setwise convergence.
\end{thm}

\textbf{Proof.} It suffices to provide a counterexample. We build on an example from \cite{YukselOptimizationofChannels} (used in a different context) in the following. Let $\mathbb{Y}=[0,1]$, $\mathbb{U}^1=\mathbb{U}^2=\{0,1\}$, and $y \sim m$ where $m$ the Lebesgue measure (uniform
distribution) on $[0,1]$. Let
\begin{eqnarray}\label{star2}
  L_{nk}= \left[\frac{2k-2}{2n},\frac{2k-1}{2n}\right), \quad
    R_{nk}= \left[\frac{2k-1}{2n},\frac{2k}{2n}\right)
\end{eqnarray}
and define the {\it square wave} function
\[
  h_n(t) = \sum_{k=1}^n\bigl(1_{\{t \in L_{nk}\}} - 1_{\{t \in R_{nk}\}} \bigr).
\]
Define further $f_n(t)=h_n(t)+1$. Let
  $B_{n,1} = \bigcup_{k=1}^n L_{nk}$ and $B_{n,2}= [0,1]\setminus
  B_{n,1}$. Define $\{Q_n\}$ as the sequence of $2$-cell quantizers
  given by
\[
Q_n(1|y)=1_{\{y \in B_{n,1}\}}, \quad Q_n(2|y)=1_{\{y\in  B_{n,2}\}} .
\]
Let \[P_n(u^1 = 1 | y) = P_n(u^2 = 1 | y) = Q_n(1|y).\]

Define $P \in {\cal P}(\mathbb{U}^1 \times \mathbb{Y} \times \mathbb{U}^2)$ as $P(a,A,b) = 1_{\{a=b\}} \frac{1}{2}m(A)$, where $a,b \in \{0,1\}$ and  $A\in \mathcal{B}([0,1])$. 

By the proof of the Riemann-Lebesgue lemma (\cite{WhZy77}, Thm.\ 12.21), observe that for all $A\in \mathcal{B}([0,1])$,
\[
\lim_{n \to \infty} \int_{A} Q_n(1|y)m(dy) = \lim_{n\to \infty} \int_0^1 \frac{1}{2}
f_n(t) \, dt = \frac{1}{2}m(A),
\]
and thus for all $A\in \mathcal{B}([0,1])$
\begin{eqnarray}
&&\lim_{n \to \infty} P_n(u^1=1,y \in A,u^2=1) \nonumber \\
&&=\lim_{n \to \infty} \int_{A}   P_n(u^1=1|y)P_n(u^2=1|y)m(dy)  \nonumber \\
&&= \lim_{n \to \infty} \int_{A}   P_n(u^1=1|y)m(dy) \nonumber \\
&& = \frac{1}{2}m(A) \nonumber \\
&&= P(1,A,1)
\end{eqnarray}
A similar property applies for $(u^1,u^2)=(0,0), (0,1)$ and $(1,0)$ so that 
\[\lim_{n \to \infty} P_n(u^1=a,y \in A,u^2=b) \to P(a,A,b) = 1_{\{a=b\}} \frac{1}{2}m(A) \]
Thus, $P_n \to P$ setwise. But even though $P_n$ satisfies the conditional independence property that $P_n(u^1=1|y,u^2) = Q_n(1|y)$, $P$ does not satisfy the conditional independence property of $u^1$ and $u^2$ given $y$: Under $P$, $u^1$ and $y$ are independent but $u^1=u^2$ and thus $P(u^1=a|y,u^2=b) = 1_{\{a=b\}} \neq {1 \over 2} = P(u^1=a|y)$. Thus, setwise (and hence weak) convergence does not preserve the conditional independence property. \hfill $\diamond$

We refer the reader to \cite{barbie2014topology}, and the references therein, for further related results. 

Some sufficient conditions for compactness of $L_R$ under the weak convergence topology are given in \cite{gupta2014existence}:
\begin{thm}\label{SuffCon1} 
Consider a static team where the action sets $\mathbb{U}^i, i \in {\cal N}$ are compact. Furthermore, if the measurements satisfy
\[P(d{\bf y}|\omega_0) = \prod_{i=1}^n Q^i(dy^i|\omega_0),\]
where $Q^i(dy^i|\omega_0) = \eta^i(y^i,\omega_0) \nu^i(dy^i)$ for some positive measure $\nu^i$ and continuous $\eta^i$ so that for every $\epsilon >0$, there exists $\delta > 0$ so that for $\rho_i(a,b)<\delta$ (where $\rho_i$ is a metric on $\mathbb{Y}^i$)
\[|\eta^i(b,\omega_0) - \eta^i(a,\omega_0) | \leq \epsilon h^i(a,\omega_0),\]
with $\sup_{\omega_0} \int h^i(a,\omega_0) \nu^i(dy^i) < \infty$, then the set $L_R(\mu)$
 is weakly compact and if $c(\omega_0,{\bf u})$ is lower semi-continuous, there exists an optimal team policy (which is deterministic and hence in $L_A(\mu)$).
  \end{thm}

The results in \cite{gupta2014existence} also apply to static reductions for sequential dynamic teams, and a class of teams with unbounded cost functions and non-compact action spaces that however satisfies some moment-type cost functions leading to a tightness condition on the set of strategic measures leading to a finite cost. In particular, the existence result applies to the celebrated counterexample of Witsenhausen \cite{wit68}. 

We also note that for a class of sequential teams with perfect recall, one can establish the existence of optimal team policies through Theorem \ref{StrategicCharacterization} and Theorem \ref{FeinbergStrategic2}. Note that the cost function is given by $c(\omega_0, {\bf u})$, where $\omega_0$ is an exogenous random variable. 
\begin{thm}\label{SuffCon2}
Consider a sequential team with a classical information structure with the further property that $\sigma(\omega_0) \subset \sigma(y^1)$. Suppose further that $\prod_k (\mathbb{Y}^k \times \mathbb{U}^k)$ is compact. If $c$ is lower semi-continuous and each of the kernels $p_n$ (defined in (\ref{kernelDefn})) is weakly continuous so that
\begin{eqnarray}\label{weakConKer}
\int f(y^n) p_n(dy^n | \omega_0,u^{1},\cdots,u^{n-1})
\end{eqnarray}
is continuous in $\omega_0,u^{1},\cdots,u^{n-1}$ for every continuous and bounded $f$, there exists an optimal team policy which is deterministic.
\end{thm}
\\
\textbf{Proof.} 
We note that when $\prod_k \mathbb{Y}^k \times \mathbb{U}^k$ is compact, the set of all probability measures with a fixed marginal on $\omega_0$ would form a weakly compact set. Therefore, it suffices to ensure the closedness of the set of strategic measures, which leads to the compactness of the set. To facilitate such a compactness condition, we first {\it expand} the information structure so that DM $k$ has access to all the previous actions $u^{1},\cdots,u^{k-1}$ as well. Later on, we will see that this expansion is redundant. With this expansion, any weakly converging sequence of strategic measures will continue satisfying (\ref{convD2}) in the limit due to the fact that there is no conditional independence property in the sequence since all the information is available at DM $k$. On the other hand, to show that for any weakly converging sequence of strategic measures satisfying (\ref{convD1}) so does the limit, it suffices to show that
\begin{eqnarray}\label{serfozRes}
\int g(y^n,\omega_0,u^{1},\cdots,u^{n-1}) p_n(dy^n | \omega_0,u^{1},\cdots,u^{n-1})
\end{eqnarray}
is continuous in $\omega_0,u^{1},\cdots,u^{n-1}$. With $f_{\omega_0,u^{1},\cdots,u^{n-1}}(y^n) := g(y^n,\omega_0,u^{1},\cdots,u^{n-1})$, it follows that as a sequence $(\omega_0,u^{1},\cdots,u^{n-1})_k \to (\omega_0,u^{1},\cdots,u^{n-1})$, and $y^n_k \to y^n$ the following holds: \[f_{(\omega_0,u^{1},\cdots,u^{n-1})_k}(y^n_k) \to f_{\omega_0,u^{1},\cdots,u^{n-1}}(y^n).\] In other words, $f_{\omega_0,u^{1},\cdots,u^{n-1}}(y^n)$ converges {\it continuously} as it is defined in \cite{serfozo1982convergence}. Thus, the continuity of (\ref{weakConKer}) ensures that (\ref{serfozRes}) holds by a generalized convergence theorem given in Theorem 3.5 of \cite{serfozo1982convergence}. Thus, the properties (\ref{convD1})-(\ref{convD2}) are preserved under weak convergence. As a result, the existence follows from Theorem \ref{existenceT}. Now, we know that an optimal policy will be deterministic as a consequence of Theorem \ref{FeinbergStrategic2}. Thus, a deterministic policy may not make use of randomization, therefore DM $k$ having access to $\{y^{k},y^{k-1},y^{k-2}, \cdots, \omega_0\}$ is informationally equivalent to him having access to $\{y^{k},(y^{k-1},u^{k-1}),(y^{k-2},u^{k-2}), \cdots, \omega_0\}$ for an optimal policy. Thus, an optimal team policy exists.

\qed

A further existence result along similar lines, for a class of static teams, is presented in Theorem \ref{SuffCon2''}

\subsection{Measurability properties of sets of strategic measures}
We have the following result.
\begin{thm}\label{FeinbergStrategicMeasurability3}
Consider a sequential (static or dynamic) team. 
\begin{itemize}
\item (i) The set of strategic measures $L_R(\mu)$ is Borel when viewed as a subset of the space of probability measures on $\Omega_0 \times \prod_k (\mathbb{Y}^k \times \mathbb{U}^k)$ under weak convergence.
\item (ii) The set of strategic measures $L_A(\mu)$ is Borel when viewed as a subset of the space of probability measures on $\Omega_0 \times \prod_k (\mathbb{Y}^k \times \mathbb{U}^k)$ under weak convergence.
\end{itemize}
\end{thm}

\textbf{Proof.}
Recall the following supporting results. The following appears in \cite{InfiniteDimensionalAnalysis} (see Theorem 15.13 in \cite{InfiniteDimensionalAnalysis} or p. 215 in \cite{Bogachev})
\begin{fact}\label{MonotoneClassRamon}
Let $\mathbb{S}$ be a Polish space and  $M$ be the set of all measurable and bounded functions $f: \mathbb{S} \to \mathbb{R}$. Then, for any $f\in M$,  the integral \[\int \pi(dx) f(x)\] defines a measurable function on ${\cal P}(\mathbb{S})$ under the topology of weak convergence.
\end{fact}

This is a useful result since it allows us to define measurable functions in integral forms on the space of probability measures when we work with the topology of weak convergence. The following result follows from Fact \ref{MonotoneClassRamon} and Theorem 2.1 of Dubins and Freedman \cite{DubinsFreedman} and Proposition 7.25 in Bertsekas and Shreve \cite{BertsekasShreve}.
\begin{fact}\label{DubinsFreedmanTheorem}
Let $\mathbb{S}$ be a Polish space. A function $F: {\cal P}(\mathbb{S}) \to {\cal P}(\mathbb{S})$ is measurable on ${\cal B}({\cal P}(\mathbb{S}))$ (under weak convergence), if for all $B \in {\cal B}(\mathbb{S})$
$(F(\cdot))(B): {\cal P}(\mathbb{S}) \to \mathbb{R}$ is measurable under weak convergence on ${\cal P}(\mathbb{S})$, that is for every $B \in {\cal B}(\mathbb{S})$, $(F(\pi))(B)$ is a measurable function when viewed as a function from ${\cal P}(\mathbb{S})$ to $\mathbb{R}$.
\end{fact}

We note that the topology considered in \cite{DubinsFreedman} is not the weak convergence topology, but since the topology considered there is not stronger than weak convergence, the result applies in this case as well.
\begin{itemize}
\item(i) {\bf Measurability of $L_R(\mu)$}
%
If $\prod_k \mathbb{Y}^k \times \mathbb{U}^k$ were compact, one could use the characterization in Theorem \ref{StrategicCharacterization}: For every given $g$, both sides of (\ref{convD1})-(\ref{convD2}) define measurable functions on ${\cal P}(\Omega_0 \times \prod_k \mathbb{Y}^k \times \mathbb{U}^k)$. One can construct a countable collection of {\it weak-convergence determining} functions on $\Omega_0 \times \prod_k \mathbb{Y}^k \times \mathbb{U}^k$, by Theorem 3.4.5 in \cite{ethier2009markov}: By Facts \ref{MonotoneClassRamon} and \ref{DubinsFreedmanTheorem}, for each continuous function the set of probability measures which satisfy the equalities in  (\ref{convD1})-(\ref{convD2}) is measurable and one only needs a countable number of such equalities, the intersection of the sets of probability measures satisfying all of these has to be Borel. If $\prod_k \mathbb{Y}^k \times \mathbb{U}^k$ is not compact, by the Borel isomorphism theorem (Appendix 1 in \cite{dynkin1979controlled} or  {\sl Chapter 13} in \cite{Dud02}), one could construct a bijection between $[0,1]$ and $\Omega_0 \times \prod_k \mathbb{Y}^k \times \mathbb{U}^k$, and one could follow the proof of Theorem 1 in \cite{dynkin1979controlled} so that only polynomial functions on $[0,1]$ are used as the test functions on this new space. By the same reasoning as the compact case, the intersection of the sets of probability measures satisfying countably many equalities of the form (\ref{convD1})-(\ref{convD2}) has to be Borel.

\item(ii) {\bf Measurability of $L_A(\mu)$}: To establish the measurability of $L_A(\mu)$, define $L^{classical}(\mu)$ as follows. Enlarge the information structure so that the information structure is classical such that the information at DM $k$ is ${\cal I}^k = \{y^k\} \bigcup \cup_{1 \leq s \leq k-1} {\cal I}^s$ for $k > 1$, and write
\begin{eqnarray}
L^{classical}(\mu) &&:= \bigg\{P \in {\cal P}\bigg(\Omega_0 \times \prod_{k=1}^N (\mathbb{Y}^k \times \mathbb{U}^k)\bigg): \nonumber \\
&& \quad \quad P(B) = \int \mu(d\omega, d{\bf y}) \prod_k 1_{\{ \gamma^k({\cal I}^k) \in B^k \}} \bigg\}
\end{eqnarray}
It follows that $L^{classical}(\mu)$ is Borel by Theorem 3.2 in Feinberg \cite{feinberg1996measurability} (which attributes the result to Blackwell \cite{blackwell1976stochastic}) since this information structure may be viewed as a fully observed Markov decision process with $x_k = {\cal I}^k$ . The measurability of $L_A(\mu)$ follows from the fact that $L_A(\mu) =L_R(\mu)  \cap L^{classical}(\mu)$, with both intersected sets being measurable.
\end{itemize}
 \qed

\begin{remark} [Implications of measurability properties and Witsenhausen's Standard Form] The Borel measurability properties of $L_R(\mu)$ is useful in establishing the {\it universal measurability} and the {\it lower semi-analytic} properties of the value functions: Define
\[J^*(\mu) =  \inf_{P \in L_R(\mu)} \int P(ds) c(s)\]
Building on the Borel measurability result for $L_R(\mu)$, and following steps similar to Dynkin and Yushkevich \cite{dynkin1979controlled} and Lemma 4.1 of Feinberg \cite{feinberg1996measurability}
 for single decision maker problems, the universal measurability of $J^*$ can be established. Such a notion is useful when Borel measurability of $J^*$ cannot be established, yet algorithms such as dynamic programming can be carried out through the verification of universal measurability properties: It suffices for a function to be universally measurable (and not necessarily Borel measurable) for its integration with respect to some probability measure to be well-defined. We refer the reader to {\sl Chapter 7} in Bertsekas and Shreve \cite{BertsekasShreve} for a comprehensive discussion on such measurability properties. In particular, according to another model for sequential teams, known as {\it Witsenhausen's Standard Form} \cite{WitsenStandard}, for optimization of sequential dynamic teams with finite horizons, a dynamic programming principle can be applied which essentially expresses the optimization problem as a terminal-stage cost function. Here, every DM acts given the policies of the previous DMs optimally. The universal measurability and the stronger condition of having the {\it lower semi-analytic} property allow for the dynamic programming recursions to be well-defined: A measurable image of a Borel set is called an {\it analytic} set and a function $f$ is called lower semi-analytic if $\{x: f(x) \leq c\}$ is analytic for each scalar $c$; see \cite{dynkin1979controlled}. It follows then that the dynamic programming recursions can be defined under Witsenhausen's Standard Form for a large class of problems. We also note that when the cost function is stage-wise additive and further assumptions are placed on the primitive variables, the analysis reduces to the usual dynamic programming formulation for state-space models. 


\end{remark}


\section{Convexity of Team Problems and Information Structures}

\subsection{Convexity of static team problems and an equivalent representation of cost functions}\label{convexStat}

We begin this section with the following definition.
\begin{defn}
A (static or dynamic) team problem is convex on ${\bf \Gamma}$ if $J(\underline{\gamma}) < \infty$ for all
$\underline{\gamma} \in {\bf \Gamma}$ and for any $\alpha \in (0,1), \underline{\gamma}_1, \underline{\gamma}_2 \in {\bf \Gamma}$:
\[J(\alpha \underline{\gamma}_1 + (1 - \alpha) \underline{\gamma}_2 ) \leq \alpha J(\underline{\gamma}_1)  + (1 - \alpha) J(\underline{\gamma}_2) \]
\end{defn}



We state the following immediate result without proof, more general refinements will be stated later in the paper. 
\begin{thm} \label{sufficientConv} 
Consider a static team. $J(\underline{\gamma})$ is convex if $c(\omega_0,\bf{u})$ is convex in ${\bf u}$ for all $\omega_0$, provided that $J(\underline{\gamma}) < \infty$ for all $\underline{\gamma} \in {\bf \Gamma}$.
\end{thm}

The condition in Theorem \ref{sufficientConv} is not tight, however, due to information structure and measurability aspects. 
\begin{example}\label{Ornek1}
Consider $\Omega= [0,1]$ and let $P$ be the uniform distribution on $\Omega$, with $N=2$, $\mathbb{U}^1 = \mathbb{U}^2 = [1,2]$. Let: 
\[c(\omega,u^1,u^2) = 1_{\{\omega \in [0,0.9]\}} \bigg((u^1 -2)^2 + (u^2-2)^2\bigg) + 1_{\{\omega \in (0.9,1]\}} \bigg(\sqrt{1+u^1} + \sqrt{1+u^2} \bigg)\]
Now, suppose further that $I^1= I^2 = \eta^1(\omega) = \eta^2(\omega) = 1_{\{\omega \in [0,0.1)\}}$. It follows that here the team problem is convex, even though $c(\omega,u^1,u^2)$ is not convex on $\{\omega: \omega \in (0.9,1]\}$, which has a non-zero probability measure. To see this, note that one may view this optimization problem as $J(u^1_1,u^1_2; u^2_1, u^2_2)$ where $u^i_j = \gamma^i(\omega_j)$, with $\omega_1 \equiv \{\omega: \omega \in [0,0.1)\}$ and $\omega_2 \equiv \{\omega: \omega \in [0.1,1]\}$. It follows that
\[J(u^1_1,u^1_2; u^2_1, u^2_2) = \sum_{i=1,2} 0.1 (u^i_1-2)^2 + 0.8 (u^i_2-2)^2 + 0.1 (\sqrt{u^i_2+1})\] 
The Hessian of $J$ is a diagonal matrix with strictly positive entries, leading to the convexity of the problem.
\end{example}

In the following, we will make use of the fact that $u^k \leftrightarrow y^k \leftrightarrow \{{\bf y}^{-k}, \omega\}$ form a Markov chain almost surely. Before proceeding further, let us note that the {\it join} of two $\sigma$-fields over some set $\mathbb{X}$ is the coarsest $\sigma$-field containing both. The {\it meet} of two $\sigma$-fields is the finest $\sigma$-field which is a subset of both. Let ${\cal F}^i$ be the $\sigma$-field generated by $\eta^i$ over $\Omega$, and let ${\cal F}_c = \bigcap_k {\cal F}^k$ be the meet of these fields, this is termed as {\it common knowledge} by Aumann \cite{Aumann} for finite probabilities spaces. In addition, let ${\cal F}_j$ be the {\it join} of the $\sigma$-field, denoted with ${\cal F}_j = \bigcup_k {\cal F}^k$.

\noindent {\bf An equivalent representation of the cost through iterated expectations.} Let us express the expected cost under a given measurable team policy $\underline{\gamma}$ as follows. We obtain from the law of the iterated expectations that
\begin{eqnarray}
&& E[c(\omega_0,{\bf u})] = E\bigg[ E[c(\omega_0,{\bf u}) | {\bf y}, {\bf u}] \bigg] = \bigg(\int P(d\omega_0 | {\bf y})  c(\omega_0, u^1,\cdots,u^N) \bigg), \label{equiRepresentation}
\end{eqnarray}
and thus with  $\tilde{c}(y^1,\cdots,y^N,u^1,\cdots,u^N) = \int P(d\omega_0 | {\bf y})  c(\omega_0, u^1,\cdots,u^N)$, the function to be minimized can be viewed as $E[\tilde{c}(y^1,\cdots,y^N,u^1,\cdots,u^N)]$. To simplify the notation, we will use the expression $E[c(\omega_0,{\bf u}) | {\bf y}]$ instead of $E[c(\omega_0,{\bf u}) | {\bf y}, {\bf u}]$ by viewing ${\bf u}$ as a fixed team control action (and not as a random variable).

\begin{thm}\label{Extension2}
\begin{itemize}
\item[(i)]If a team problem is convex, then $E[c(\omega_0, {\bf u}) | {\cal F}_c]$ is convex in ${\bf u}$ almost surely.  \\
\item[(ii)] If $E[c(\omega_0, {\bf u}) | {\cal F}_j]$ is convex in ${\bf u}$ almost surely, then the team problem is convex on the set of team policies that satisfy $J(\underline{\gamma}) < \infty$.
\end{itemize}
\end{thm}

\textbf{Proof.}
(i) We will show the contra-positive. Let $B$ be a Borel set such that $P(B) > 0$, $B \in {\cal F}_c$, and $E[c(\omega_0, {\bf u}) | B]$  be non-convex so that there exist ${\bf u}$ and ${\bf u}'$ and $\lambda \in (0,1)$ such that
\[ E[c(\omega_0, \lambda {\bf u} + (1 - \lambda){\bf u}') | B] >  \lambda E[c(\omega_0,  {\bf u}) | B]  +  (1 - \lambda) E[c(\omega_0, {\bf u}') | B] \]
Now, let $\underline{\gamma}$ and ${\bf \gamma}$ be two team policies so that these only differ on $B$; and on $B$ $\underline{\gamma}={\bf u}$ and ${\bf \gamma} = {\bf u}'$. Such measurable policies exist, for example by taking $\underline{\gamma}(\omega) = \{0,0,\cdots,0\}$ when $\omega \notin B$. These policies are both Borel measurable and are admissible given the information structure. Then $J(\lambda \underline{\gamma} + (1 - \lambda){\bf \gamma} ) > \lambda J( \underline{\gamma}) + (1 - \lambda) J( {\bf \gamma})$. \\

(ii) 
We adopt the equivalent representation (\ref{equiRepresentation}) in this part of the proof. Note that under any measurable policy, the random variable $\tilde{c}(y^1,\cdots,y^N,u^1,\cdots,u^N)$ is measurable on the $\sigma$-field generated by ${\bf y}$ and thus the join $\sigma$-field. The proof then follows from the following. Consider two policies $\underline{\gamma}$ and ${\bar{\underline{\gamma}}}$ with finite expected costs. It follows then that
\begin{eqnarray*}
&&J(\lambda \underline{\gamma} + (1 - \lambda) {\bar{\underline{\gamma}}})  \\
&&= \int P(d{\bf y}) \tilde{c}(y^1,\cdots,y^N, \lambda \gamma^1(y^1) + (1 - \lambda) \bar{\gamma}^1(y^1),\cdots,\lambda \gamma^N(y^N) + (1 - \lambda) \bar{\gamma}^N(y^N)) \\
&& \leq  \int P(d{\bf y}) \bigg(\lambda \tilde{c}(y^1,\cdots,y^N, \gamma^1(y^1) ,\cdots, \gamma^N(y^N) ) \\
&& \quad \quad \quad \quad + (1-\lambda) \tilde{c}(y^1,\cdots,y^N, \bar{\gamma}^1(y^1),\cdots, \bar{\gamma}^N(y^N) ) \bigg) \\
&&= \lambda  J(\underline{\gamma}) + (1 - \lambda) J({\bar{\underline{\gamma}}}) 
\end{eqnarray*}

\qed

It can be observed that Example \ref{Ornek1} satisfies the conditions of Theorem \ref{Extension2}. We will use these to study Witsenhausen's counterexample \cite{wit68} later in the paper. 

\subsubsection{A generalization of Radner and Krainak et. al.'s theorems}

We provide a generalization of Radner's or Krainak et al.'s theorem by utilizing an information structure dependent nature of convexity. For example, Radner or Krainak et.al's theorems are not applicable to Example \ref{Ornek1}. 


\begin{thm}\label{Thm:Y4}
Let $\{ J; \Gamma^i , i\in\IN \}$ be a static stochastic team
problem, the loss function $E[c( \omega_0 , \boldu ) | {\cal F}_j]$ is convex and continuously differentiable
in $\boldu$ almost surely. Let $\undergamma^* \in \boldGamma$ be a policy $N$-tuple with a
finite cost $(J(\undergamma^*) < \infty)$, and suppose that for
every $\undergamma \in \boldGamma$ such that $J(\undergamma) <
\infty$, the following holds: 
\begin{equation} 
\sum_{i \in {\cal N}} E \{ \nabla_{u^i} \tilde{c}({\bf y}, \undergamma^* (\boldy)) [ \gamma^i (y^i)
- \gamma^{i*} (y^i)] \} \geq 0, \label{eq:20}
\end{equation}
where $\tilde{c}({\bf y},{\bf u}) := E[c( \omega_0 , \boldu ) | {\cal F}_j]$. Then,
$\undergamma^*$ is a {\it team-optimal policy}, and it is {\it
unique} if $\tilde{c}({\bf y},{\bf u})$ is strictly convex in $\boldu$ almost surely.
\end{thm}

\noindent
{\bf Proof.} The proof follows by defining the new loss function $\tilde{c}({\bf y},{\bf u}) = E[c( \omega_0 , \boldu ) | {\cal F}_j]$. The result then follows as in Theorem \ref{Thm:TB4}. \qed

\begin{thm}\label{Thm:Y5} Let $\{ J; \Gamma^i, i\in\IN \}$ be a static
stochastic team problem which satisfies all the hypotheses of
Theorem~\ref{Thm:Y4}, with the exception of inequality (\ref{eq:20}).  Instead
of~(\ref{eq:20}), let either (c.5) or (c.6) be satisfied with $c$ replaced with $\tilde{c}$.  Then, if
$\undergamma^* \in \boldGamma$ is a stationary policy it is also
team optimal.  Such a policy is unique if $E[c(\omega_0, \boldu ) | {\cal F}_j]$ is
strictly convex in $\boldu$, a.s.
\end{thm}

\noindent
{\bf Proof.} The proof follows by defining the new loss function $\tilde{c}$ as in the proof of Theorem \ref{Thm:Y4}, and following Theorem  \ref{Thm:TB5}. \hfill $\diamond$

\subsection{A further existence result for static teams through the equivalent representation}\label{furtherExis}

By defining $\tilde{c}$ as the new loss function, we can obtain an existence result which allows for Theorem \ref{SuffCon2} to be applicable to a static team where the information is nested (and thus the information structure is classical). These relax some of the conditions presented in \cite{AbhishekPhD}.

\begin{thm}\label{SuffCon2''}
Consider a static team with a classical information structure (that is, with an expanding information structure so that $\sigma(y^n) \subset \sigma(y^{n+1}), n \geq 1$). Suppose further that $\prod_k (\mathbb{Y}^k \times \mathbb{U}^k)$ is compact. If 
\[\tilde{c}(y^1,\cdots,y^N,u^1,\cdots,u^N):=E[c(\omega_0,{\bf u}) | {\bf y}, {\bf u}]\]
is jointly lower semi-continuous in ${\bf u}$ for every ${\bf y}$,
and for every $n > 1$,
\begin{eqnarray}\label{weakConKer''}
\int f(y^n) P(dy^n | y^{n-1})
\end{eqnarray}
is continuous in $y^{n-1}$ for every continuous and bounded $f$, there exists an optimal team policy which is deterministic.
\end{thm}
\\
\textbf{Proof.} 
Different from Theorem \ref{SuffCon2}, we eliminate the use of $\omega_0$, and study the properties of the set of strategic measures. Also, instead of the weak topology, we will use the $w$-$s$ topology introduced by Sch\"al \cite{Schal}: The $w$-$s$ topology on the set of probability measures ${\cal P}(\mathbb{X} \times \mathbb{U})$ is the coarsest topology under which $\int f(x,u) \nu(dx,du): {\cal P}(\mathbb{X} \times \mathbb{U}) \to \mathbb{R}$ is continuous for every measurable and bounded $f$ which is continuous in $u$ for every $x$ (but unlike weak topology, $f$ does not need to be continuous in $x$).

As in the proof of Theorem \ref{SuffCon2}, when $\prod_k \mathbb{Y}^k \times \mathbb{U}^k$ is compact, the set of all probability measures on $\prod_k \mathbb{Y}^k \times \mathbb{U}^k$ forms a weakly compact set. Since the marginals on $\prod_k \mathbb{Y}^k$ is fixed, \cite[Theorem 3.10]{Schal} (or  \cite[Theorem 2.5]{balder2001}) establishes that the set of all probability measures with a fixed marginal on $\prod_k \mathbb{Y}^k$ is relatively compact under the $w$-$s$ topology. Therefore, it suffices to ensure the closedness of the set of strategic measures, which leads to the sequential compactness of the set under this topology. To facilitate such a compactness condition, as earlier we first {\it expand} the information structure so that DM $k$ has access to all the previous actions $u^{1},\cdots,u^{k-1}$ as well. Later on, we will see that this expansion is redundant. With this expansion, any w-s converging sequence of strategic measures will continue satisfying (\ref{convD2}) in the limit due to the fact that there is no conditional independence property in the sequence since all the information is available at DM $k$. That is, $P_n(du^n | y^n,y_{[0,n-1]},u_{[0,n-1]})$ satisfies the conditional independence property trivially as all the information is available. On the other hand, for each element in the sequence of conditional probability measures, the conditional probability for the measurements writes as
$P(dy^n | y_{[0,n-1]},u_{[0,n-1]}) = P(dy^n | y^{n-1})$.  We wish to show that this also holds for the $w$-$s$ limit measure. Now, we have that for every $n$, $y^n \leftrightarrow y^{n-1} \leftrightarrow u_{[0,n-1]}$ forms a Markov chain. By considering the convergence properties only on continuous functions, as in (\ref{convD1}), if $P(y^n | y^{n-1})$ is weakly continuous, the $w$-$s$ limit (and thus the weak limit also, since weak convergence is weaker than $w$-$s$ convergence) will also satisfy this property. 

Thus, (\ref{convD1}) is also preserved by the weak continuity condition (\ref{weakConKer''}) as in (\ref{serfozRes}) in the proof of Theorem \ref{SuffCon2}. Hence, for any $w$-$s$ converging sequence of strategic measures satisfying (\ref{convD1})-(\ref{convD2}) so does the limit since the team is static and with perfect-recall. By \cite[Theorem 3.7]{Schal}, and the generalization of Portmanteau theorem for the $w$-$s$ topology, the lower semi-continuity of $\int \mu(d{\bf y}, d{\bf u}) \tilde{c}({\bf y},{\bf u})$ over the set of strategic measures leads to the existence of an optimal strategic measure. As a result, the existence follows from similar steps to that of Theorem \ref{existenceT}. Now, we know that an optimal policy will be deterministic as a consequence of Theorem \ref{FeinbergStrategic2}. Thus, a deterministic policy may not make use of randomization, therefore DM $k$ having access to $\{y^{k},y^{k-1},y^{k-2}, \cdots\}$ is informationally equivalent to him having access to $\{y^{k},(y^{k-1},u^{k-1}),(y^{k-2},u^{k-2})\}$ for an optimal policy. Thus, an optimal team policy exists.
\qed

\subsection{Convexity of Sequential Dynamic Teams}\label{convexDynam}

\subsubsection{Convexity of the reduced model}

The static reduction of a sequential dynamic team problem, if exists, is not unique. However, the following holds: Either all of the static reductions are convex or none is. This argument follows from the fact that the policy spaces are isomorphic in the sense discussed earlier in Section \ref{EquivIS}, and the induced costs are identical given a fixed team policy under the static reduction and the original formulation. The following is then immediate.
\begin{thm}
A stochastic dynamic team problem with a static reduction is convex if and only if its static reduction is.
\end{thm}



\begin{remark} Consider the condition discussed in Theorem \ref{Extension2}. If the expression $E[c(\omega_0 , {\bf u}) | {\cal F}_j]$ is convex almost surely in one reduction given by (\ref{staticReduc})-(\ref{staticReduc2}), so it is in another, since the reduced models are such that the measure induced on the measurements are mutually absolutely continuous.
\end{remark}

{\bf Non-convexity of the Witsenhausen counterexample and its variants.} Consider the celebrated Witsenhausen's counterexample \cite{wit68}: This is a dynamic non-classical team problem with $y^1$ and $w^1$ zero-mean independent Gaussian random variables with unit variance and $u^1=\gamma^1(y^1)$, $u^2 = \gamma^2(u^1+w^1)$ and the cost function $c(\omega,u^1,u^2)= k^2(y^1-u^1)^2 + (u^1-u^2)^2$ for some $k > 0$: The static reduction proceeds as follows:
\begin{eqnarray}\label{stW}
&&\int (k(u^1-y^1)^2 + (u^1-u^2)^2) Q(dy^1)\gamma^1(du^1|y^1)\gamma^2(du^2|y^2)P(dy^2|u^1) \nonumber \\
&& = \int \bigg( (k(u^1-y^1)^2 + (u^1-u^2)^2) \gamma^1(du^1|y^1)\gamma^2(du^2|y^2) {\eta(y^2-u^1) dy^2 \over \eta(y^2)} \bigg) Q(dy^1) \eta(y^2)dy^2 \nonumber \\
&& = \int \bigg( (k(u^1-y^1)^2 + (u^1-u^2)^2) \gamma^1(du^1|y^1)\gamma^2(du^2|y^2) {\eta(y^2-u^1) dy^2 \over \eta(y^2)} \bigg) Q(dy^1) Q(dy^2) \nonumber \\
\end{eqnarray}
where $Q$ denotes a zero-mean Gaussian measure with unit variance, and $\eta$ its density. Another interesting example is the point-to-point communication problem: Here, the setup is exactly as in the Witsenhausen's counterexample, but $c(\omega,u^1,u^2)= k^2(u^1)^2 + (y^1-u^2)^2$. This problem is a peculiar one in that, even though the information structure is non-classical, and is non-convex; an optimal encoder and decoder is linear. A proof of this result builds on information theoretic ideas, such as the data-processing inequality (see {\sl Chapters 3, 11} in \cite{YukselBasarBook} for a detailed account). In this case, the reduction writes as:
\begin{eqnarray}\label{stQ}
&&\int (k(u^1)^2 + (y^1-u^2)^2) Q(dy^1)\gamma^1(du^1|y^1)\gamma^2(du^2|y^2)P(dy^2|u^1) \nonumber \\
&& = \int \bigg(  (k(u^1)^2 + (y^1-u^2)^2) \gamma^1(du^1|y^1)\gamma^2(du^2|y^2) {\eta(y^2-u^1) dy^2 \over \eta(y^2)} \bigg) Q(dy^1) Q(dy^2) \nonumber \\
\end{eqnarray}


Consider the static reduction of Witsenhausen's counterexample and the Gaussian signaling problem (\ref{stW})-(\ref{stQ}). 
For both (\ref{stW}) and (\ref{stQ}), using the fact that $e^{-x^2}$ is not a convex function, we recognize that this problem is not convex by Theorem \ref{Extension2}(i) (with the common knowledge/information being the trivial $\sigma$-algebra consisting of the empty set and its complement). 

We note that Witsenhausen states without proof in \cite{wit68} (p. 134) that the counterexample is non-convex in $\gamma^1$ for every optimal $\gamma^2$ (selected as a best response to $\gamma^1$). The discussion above can be viewed as an explicit proof for this result. Note also that for both problems above, linear policies contain person by person optimal policies, but this does not imply global optimality. For the first problem, Witsenhausen had shown the suboptimality of linear policies. For the second problem (\ref{stQ}), however, linear policies are indeed optimal. 

%

%
%
%
%

\subsubsection{Partially nested information structures: Convexity of the reduced model}\label{QuasiRed}

As reviewed earlier in Section \ref{witsenInfoStructureReview}, an important information structure which is not nonclassical, is of the {\it partially nested} type. For such team problems with partially nested information, a static reduction exists under certain invertibility conditions as discussed earlier. For such problems, the cost function is not altered by the static reduction. This leads to the following result. 



\begin{thm}
Consider a partially nested stochastic dynamic team which admits a static reduction where the cost function $c(\omega_0 , {\bf u})$ convex in ${\bf u}$. The team problem is convex.
\end{thm}

We note that Ho and Chu \cite{HoChu} established this result that for the special setup involving the partially nested LQG teams. In this case, optimal policies are linear through an equivalence to static teams.

\subsubsection{Stochastic partial nestedness: A probabilistic definition of nestedness, its relation to convexity and signaling}

When the information structure is non-classical, the decision makers may use their actions to communicate with each other. This phenomenon is known as signalling. When signaling is present, the problem has a communications flavour and any communication problem is inherently non-convex, see Theorem 4.3.1 in \cite{YukselBasarBook}. It is known that quasi-classical information structures eliminate the incentive for signaling, since the future decision makers already have access to the information at the signaling decision maker (see \cite{Rantzer}, but also \cite{CasalinoDavoli}, \cite{Voulgaris}, \cite{Bamieh} and \cite{Rotkowitz} among other papers). On the other hand, one can also put a probabilistic flavor: \cite{YukTAC09} identified such a probabilistic, but rather restrictive, characterization; see also \cite{MahajanYuksel2010}. In the following, we exhibit that the static reduction provides an effective method to identify when lack of a signaling incentive can be established and can lead to a more refined probability and information structure dependent characterization of {\it nestedness}, that encompasses partial nestedness which is a probability-free characterization. 

\begin{defn}
The information structure of a sequential team problem is {\it stochastically partially nested}, if for an arbitrary cost function $c: \Omega_0 \times \prod_k \mathbb{U}^k \to \mathbb{R}$ there exists a static reduction of this team which does not alter the loss function.  
\end{defn}

This definition implies the following result.
\begin{lem}
Consider a sequential team problem with a stochastically partially nested information structure. If the cost function $c(\omega_0 , {\bf u})$ is convex in ${\bf u}$, then the team problem is convex.
\end{lem}

\textbf{Proof.} The static reduction of this team preserves convexity of the loss function, for an arbitrary convex loss function $c: \Omega_0 \times \prod_k \mathbb{U}^k \to \mathbb{R}$. Thus, the reduced problem, and hence the original problem is convex. 
 \qed

This definition essentially requires that there is no active information transmission and there is no signaling incentive. To provide an example of an information structure that is not partially nested, but that preserves convexity properties and is stochastically partially nested, let us consider the following example from \cite{YukTAC09}, with the approach of this paper. 

\begin{defn}\label{decoupled}
Let ${\bf x}=\begin{bmatrix}x^1 & x^2 & \dots & x^L\end{bmatrix}^T,$ 
and a system evolve as
\[x^i_{t+1}= f^i(x^i_t,u^i_t,w^i_t), \quad \quad \quad y^i_t = g^i({\bf x}_t, v^i_t), \quad i \in \{1,2,\dots,L\},\]
\[z_{t+1}= f(z_t,u^1_t,u^2_t,\dots,u^L_t),\]
for measurable functions $f^i, g^i, f$ and $\{w^i_t\}$ independent state disturbance processes and $\{v^i_t\}$, observation noise processes for $i \in \{1,2,\dots,L\}$. In the above $(x^i_t,w^i_t,v^i_t, i \in \{1,2,\dots,L\},z_0)$ are independent second-order processes. Suppose each of the DMs has access to the additional $\{z_t\}$ process: $\tilde{y}^i_t=\{y^i_t , z_t\}, \quad i \in \{1,2,\dots,L\}$. If the information available at each controller is such that
\[x^i_t \leftrightarrow (y^i_{[0,t]},u^i_{[0,t-1]}) \leftrightarrow \{x^j_0,z_0, w^j_{[0,t-1]},y^{j}_{[0,t]}, \quad j \neq i \}\]
form Markov chains, so that almost surely,
\[P\bigg(dx^i_t | y^i_{[0,t]},u^i_{[0,t-1]}\bigg) = P\bigg(x^i_t | y^i_{[0,t]},u^i_{[0,t-1]}, \{x^j_0,z_0, w^j_{[0,t-1]},y^{j}_{[0,t]}, \quad j \neq i \} \bigg)\]
 for all $t$ and $i$, then such an information structure is said to be stochastically decoupled.
\end{defn}

%

Note that this information structure is non-classical since the DMs can signal information to each other, yet there is no nestedness in the information $\sigma$-fields.
\begin{thm}\label{stochDecoupled}
Let there be an optimization problem with the objective to be minimized as: $E[\sum_{t=0}^{T-1} c_1(x^1_t,u^1_t) + c_2(x^2_t,u^2_t) + \dots c_L(x^L_t,u^L_t)].$ If the controllers have stochastically decoupled information structures, the optimization problem can be cast as $L$ decoupled optimization problems with a classical information structure.
\end{thm}

\textbf{Proof:}
As earlier, through the law of the iterated expectations
\begin{eqnarray}
&& E[\sum_{t=0}^{T-1} c_1(x^1_t,u^1_t) + c_2(x^2_t,u^2_t) + \dots c_L(x^L_t,u^L_t) ] \nonumber \\
&& = E[ E[\sum_{t=0}^{T-1} c_1(x^1_t,u^1_t) + c_2(x^2_t,u^2_t) + \dots c_L(x^L_t,u^L_t)  |{\bf{\tilde{y}}}_{t},{\bf I}_{t-1}] ]. \nonumber
\end{eqnarray}
Here $I^i_k = \{\tilde{y}^i_s, s \leq k\}$ and ${\bf I}_k := \{I^1_k,\cdots, I^L_k \}$. We will show that due to the Markov chain condition in Definition \ref{decoupled}, the above writes as: $E[\sum_{t=0}^{T-1} \sum_{i=1}^L E[c_i(x^i_t,u^i_t) | y^i_t, I^i_{t-1})]$. This follows from the following observation. Let ${\bf x} = \{x^1, x^2,\dots,x^L \}$. Suppose in the following we first assume that the conditional densities exist almost surely and are denoted by $p(\cdot | \cdot)$. The conditional densities then write (almost surely) as:
\begin{eqnarray}
&&p({\bf x}_t|{\bf{\tilde{y}}}_{t},{\bf I}_{t-1}) = {p({\bf x}_t,{\bf{\tilde{y}}}_{t}|{\bf I}_{t-1}) \over \int_{\mathbb{X}} p({\bf x}_t,{\bf{\tilde{y}}}_{t}|{\bf I}_{t-1})} = {p(z_t|{\bf I}_{t-1}) p({\bf x}_t,{\bf y}_{t}|{\bf I}_{t-1}) \over \int_{\mathbb{X}} p(z_t|{\bf I}_{t-1}) p({\bf x}_t,{\bf y}_{t}|{\bf I}_{t-1}) } \nonumber \\
&&= {p(z_t|{\bf I}_{t-1}) p({\bf x}_t,{\bf y}_{t}|{\bf I}_{t-1}) \over  p(z_t|{\bf I}_{t-1}) \int_{\mathbb{X}} p({\bf x}_t,{\bf y}_{t}|{\bf I}_{t-1})} = p({\bf x}_t|{\bf y}_{t},{\bf I}_{t-1}) =\prod_{i=1}^L p(x^i_t|{\bf y}_{t},x^1_t,x^2_t, \dots,x^{i-1}_t, {\bf I}_{t-1}) \nonumber \\
%
%
&&=\prod_{i=1}^L P\bigg(dx^i_t|{\bf y}_{t},(F^1_t(x^1_0,w^1_{[0,t-1]},{\bf I}_{t-1})), \dots,(F^{i-1}_t(x^{i-1}_0,w^{i-1}_{[0,t-1]},{\bf I}_{t-1})), {\bf I}_{t-1}\bigg) \nonumber \\
&&=\prod_{i=1}^L p(x^i_t|y^i_t, I^i_{t-1}) \label{indep}
 \end{eqnarray}
If densities do not exist, the Radon-Nikodym derivative of $P(d{\bf x}_t, d{\bf{\tilde{y}}}_{t}|{\bf I}_{t-1})$ and $P(d{\bf{\tilde{y}}}_{t} | {\bf I}_{t-1})$ ensures that a conditional probability measure $P(d{\bf x}_t|{\bf{\tilde{y}}}_{t},{\bf I}_{t-1})$ can be defined so that
\[P(d{\bf x}_t, d{\bf{\tilde{y}}}_{t}|{\bf I}_{t-1}) =  P(d{\bf x}_t|{\bf{\tilde{y}}}_{t},{\bf I}_{t-1}) P(d{\bf{\tilde{y}}}_{t} | {\bf I}_{t-1}) \]
and the analysis above can be carried out by studying the properties of $P(d{\bf x}_t|{\bf{\tilde{y}}}_{t},{\bf I}_{t-1})$, and for other terms.

In the analysis above, (\ref{indep}) follows from the Markov chain hypothesis of the theorem and the rest of the arguments uses the properties of total probability. In the above, we write $x^i_t = F^i_t(x^i_0,w^i_{[0,t-1]},{\bf I}_{t-1})$, to express the explicit dependence on the variables. Since the cost function is also decoupled, it turns out that one can write the dynamic programming recursions as $L$ decoupled optimization problems. Given the new expression, $E[\sum_{t=0}^{T-1} \sum_{i=1}^L E[c_i(x^i_t,u^i_t) | y^i_t, I^i_{t-1}],$
we can write the cost as: $E[\sum_{t=0}^{T-1} \sum_{i=1}^L E[c_i(x^i_t,u^i_t)] ]$. Thus, the cost can be written as a summation of decoupled classical systems. 

 \hfill $\diamond$
\begin{example}
An example is the following. Consider a dynamical system described by
\begin{eqnarray}
&& x^1_{t+1}=a_1x^1_t + u^1_t + w^1_t, \quad x^2_{t+1}=a_2x^2_t + u^2_t + w^2_t, \quad x^3_{t+1}=a_3x^3_t + u^1_t + u^2_t + w^3_t \nonumber \\
&& y^1_{t}=(x^1_t+v^{1}_t, x^2_t + v^2_t+ v^{21}_t, x^3_t+v^{31}_t), \quad y^2_{t}=(x^1_t + v^1_t + v^{12}_t, x^2_t + v^2_t, x^3_t+v^{32}_t), \nonumber
\end{eqnarray}
where all the external random variables are zero-mean Gaussian, and the goal is the minimization of $E\bigg[\sum_{t=0}^{T-1} \bigg((x^1_t)^2 + (x^2_t)^2 + \rho_1 (u^1_t)^2 + \rho_2 (u^2_t)^2 \bigg) \bigg]$, with $\rho_1, \rho_2 >0$ constants. The control actions are measurable with the sigma-algebra generated by their causal observations and past controls: $I^i_t=\{y^i_t,I^{i}_{t-1}\}$, with $I^i_{0}={y^i_0}$. This system has a non-classical information structure, as controller 1 affects the observation at controller 2 (the third state $x^3$ acts as a communications medium between the controllers), but controller 2 cannot recover the information at controller 1. However, the optimal team policy is linear by stochastic partial nestedness and the static reduction.
\end{example}

\section{Conclusion}

In this paper, strategic measures for stochastic team problems have been introduced and properties such as convexity, and compactness and Borel measurability under weak convergence topology are studied. Sufficient conditions for each of these properties have been presented, where these conditions lead to existence of and structural results for optimal policies. It is established that team problems do not lead to a convex set of strategic measures, even in the presence of common or private independent randomness, and the global optimality of deterministic policies among the possibly randomized class has been established using such a convex analysis under an expected cost criterion. 
The problem of when a sequential team problem is convex is studied, and necessary and sufficient conditions for convexity of problems which include teams with a non-classical information structure are presented. Building on these results, generalizations of Radner's theorem have been reported. For dynamic teams, static reduction is provided as a useful method to establish not only convexity, but also obtain a systematic method through which a probabilistic relaxation for partial nestedness can be established.

\section{Acknowledgements}

We gratefully acknowledge technical discussions with Profs. Tamer Ba\c{s}ar, Abhishek Gupta and Tam\'as Linder.

\end{document}